\documentclass[10pt]{article}
\usepackage{amsmath,amsthm,epsf,graphicx}
\usepackage{graphicx}
\usepackage{amssymb}
\usepackage{epsfig}
\def\ZZ{{\mathbf Z}}
\def\NN{{\mathbf N}}

\def\hpsi{{\widehat \psi}}

\def\md{{\rm mod}\,}
\newtheorem{assumption}{Assumption}

\newtheorem{corollary}{Corollary}
\newtheorem{proposition}{Proposition}

\def\beqns{\begin{eqnarray*}}
\def\eeqns{\end{eqnarray*}}
\def\beqn{\begin{eqnarray}}
\def\eeqn{\end{eqnarray}}
\def\order{{\cal O}}

\def\C{{\bf C} \,}

                        \newcommand{\Z}{\makebox[0.04cm][l]{\sf
Z}\!{\sf Z}}
                        
                        \newcommand{\R}{{\rm I}\!{\rm R}}

\newcommand{\EQ}[2]{\begin{equation}{{#2}\label{#1}} \end{equation}}
                        \newcounter{theorem}

\renewcommand{\thetheorem}{\arabic{section}.\arabic{theorem}}

\newenvironment{thm}[2]{\begin{sloppypar}\refstepcounter{theorem}%
                        {\bf #1 \thetheorem.}\label{#2}\em{}}%
                        {\end{sloppypar}}
                        \newcommand{\theo}[3]{\begin{thm}{#1}{#2}
#3\end{thm}}

\begin{document}

\title{Uniform Convergence of Multilevel Stationary Gaussian Quasi-Interpolation}
\author{Jeremy Levesley\footnote{Department of Mathematics, University of Leicester, LE1 7RH, UK. {\tt jl1@le.ac.uk}} and Simon Hubbert\footnote{Department of Economics, Mathematics and Statisitcs, Birkbeck, University of London, WC1H 7HX, UK. {\tt s.hubbert@bbk.ac.uk}}}

\maketitle

\begin{abstract}
It is well-known that polynomial reproduction is not possible when approximating with Gaussian kernels. Quasi-interpolation schemes have been developed which use a finite number of Gaussians at different scales, which then reproduce polynomials of low degree \cite{beatson}, and thus achieve polynomial orders of convergence. At the same time, interpolation with kernels of fixed width suffers from an explosion in condition number, and information from all data points influences the approximation at any one data point (no localisation). In \cite{HL1} the authors show that, for periodic convolution with the Gaussian kernel, a multilevel scheme can give orders of approximation faster than any polynomial. In this paper we present a new multilevel quasi-interpolation algorithm, the discrete version of the algorithm in \cite{HL1}, which mimics the continuous algorithm well, to single precision accuracy, and gives excellent convergence rates for band limited periodic functions. In this paper we explain how the algorithm works, and why we achieve the numerical results we do. The estimates developed have two parts, one involving the convergence of a low degree polynomial truncation term  and  one involving  the control of the remainder of the truncation as the algorithm proceeds.
\end{abstract}

\section{Introduction}

Gaussian radial basis functions (RBFs) are, theoretically, an excellent tool for approximation in high dimensions. They are positive definite, ensuring that interpolation matrices are invertible, and theoretical errors are exponentially small as a function of the number of data points. On the other hand, in numerical calculation we meet a number of hurdles. More practically, the interpolation matrices arising from approximation using a fixed width of Gaussian become numerically singular, and efforts need to be made to stabilize the approximation process.

Motivated by numerical success in using multilevel sparse grid interpolation, quasi-interpolation and collocation methods for solving approximation and partial differential equations, we develop the analysis for multilevel quasi-interpolation of periodic functions, on uniform data, as this allows us to use Fourier techniques. In \cite{HL1} the continuous convolution case is investigated, and we conclude that for sufficiently smooth functions we achieve convergence rates faster than any polynomial.

In this paper we will construct our approximation, to a target function $f,$ via Schoenberg's approach \cite{schoen} to quasi-interpolation where, for a positive integer $n$ (which we will refer to as the sample rate) we  consider the following stationary quasi-interpolant
\EQ{quasi}{
Q_{\frac{1}{n}}(f)(x):=\sum_{\ell \in \ZZ}f\left(\frac{\ell}{n}\right)\psi\left(nx-\ell\right),\,\,\,{\rm{where}}\,\,\,\,\psi(x)=\frac{1}{\sqrt{2\pi}}\exp\left(-\frac{x^{2}}{2}\right)\,\,\,\,{\rm{(the}}\,\,{\rm{Gaussian}}\,\,{\rm{kernel).}}
}
This can be viewed as the discrete convolution, and we show that the convergence of this scheme mirrors that in \cite{HL1} well.

It is well-known that we cannot reproduce a constant when quasi-interpolating with Gaussians. In order to mitigate this problem Beatson and Light \cite{beatson} use a linear combination of Gaussians at various scales in order to satisfy the Strang-Fix conditions, and hence reproduce polynomials of a fixed degree. This gives rise to a scheme which converges at a polynomial rate. The higher the rate required, the larger the set of points, or the number of different Gaussian scales for the linear combination must be.

We choose a simpler path, and use a multilevel method.
To describe the approximation error of (\ref{quasi}) to the target function $f$ we write
\[
E_{\frac{1}{n}}(f)=Q_{\frac{1}{n}}(f)-f.
\]
We begin, at level one, by forming the quasi-interpolant (\ref{quasi}) to the target function with a starting sample rate $n.$ Then, as we move from one level to the next, we double the sample rate and form the  quasi-interpolant to the residual function (from the previous level), this is then added to the current approximation. Continuing in this  we build up an approximation to $f$ and the algorithm terminates when the residuals are sufficiently small. The residual at level one is simply $E_{\frac{1}{n}}(f)$ and the residual error at the subsequent level is defined recursively from here. Specifically, we write
\[
M_{\frac{1}{n},p}(f)=E_{\frac{1}{2^{p-1}}\frac{1}{n}}M_{\frac{1}{n},p-1}(f),
\]
to denote the multilevel error after $p$ iterations of the algorithm when  applied to $f$ with $n$ as the initial sample rate. Note that
 we set $M_{\alpha,0}=I$ to be the identity  operator so that $M_{\frac{1}{n},1}(f)=E_{\frac{1}{n}}(f).$
 We note that at any specific level we are only using information about the error locally, and so we do not require any large stencil in order to construct our approximation.

In this paper we will focus attention on target functions that belong to the cosine family $c_{m}(x)=\cos(2\pi m x),$ $m \in \ZZ$ (the same techniques can be applied to the sine family too) and our aim is describe how the magnitude of multilevel errors (i.e., $\|M_{\frac{1}{n},p}c_{m}\|_{\infty}$) behaves as the algorithm progresses. The paper is computationally intensive, so in order to explain the steps in the proof it is instructive to see how the quasi-interpolation process performs with increasing $n$.

The figure below shows a fit of the profile of $\|E_{\frac{1}{n}} c_8\|_{\infty}$ for increasing $n$. The key observation to note is that as the frequencies grow from zero (the constant function) to half of the sample rate, the approximation error grows but not beyond $1.$ For cosines whose frequency is greater than half of the sample rate the approximation errors roughly lie between $1$ and $2,$ and, due to aliasing,   the profile exhibits an oscillatory pattern as the frequencies grow.

\begin{center}
\begin{figure}[h]\caption{Profile of the upper  bound on the quasi-interpolation error to $c_{m}$ when the sample rate  $=8.$}
        \includegraphics{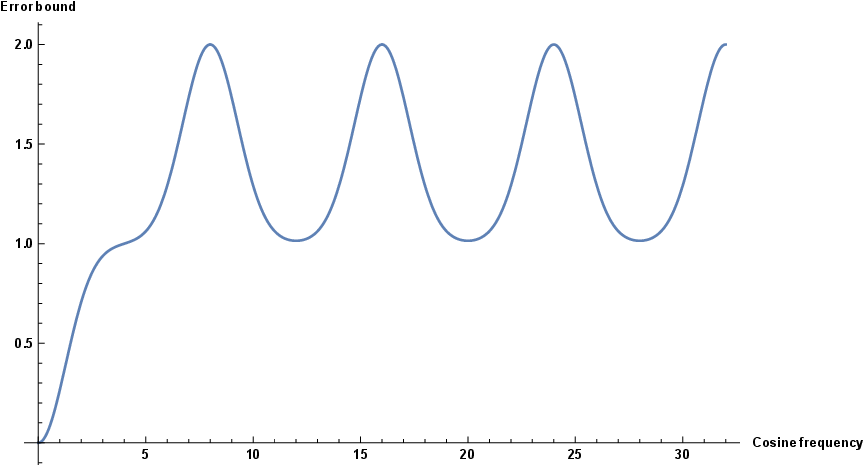}
        \label{figerr}
    \end{figure}
    \end{center}
The rationale for multilevel approximation now becomes more clear. At each level of the multilevel scheme we double the sample rate, so that more frequencies are in the left hand part of the curve, where we get reduction in error. For higher frequencies (which have smaller Fourier coefficients in the series expansion) the algorithm essentially has no effect. So we divide our analysis broadly into three parts: (a)  the cosine frequency is less than half the sample rate (less than $n=4$ in the picture above), where we get reduction in error; (b) the cosine frequency is greater than the sample rate, (greater than $n=8$ in the picture above), where the error lies between $1$ and $2,$ and  we need to keep track of the error carefully, until the point in the algorithm where the error reduces; (c) the cosine frequency is greater than half the sample rate but does not exceed the sample rate (between $4$ and $8$ in the picture)  this marks the transition region between the part where the error reduces and the part where it oscillates between $1$ and $2,$ and, as we shall see, this has a specific challenge in terms of bounding errors.

The paper is organised as follows. In Section~\ref{overview} we give a brief overview of previous work on multilevel methods, then  in Section~\ref{preliminaries} we will compose the key mathematical results that will be useful in the subsequent analysis. We begin by describing the space of periodic functions to which our algorithm will be applicable and, for a given sample rate $n$, we  establish a general bound on the quasi-interpolation error of such functions. We then pay attention to the case where the target function is a member of the cosine family $c_{m}(x)=\cos(2\pi m x),$ $m \in \ZZ.$ We derive a formula for the quasi-interpolant $Q_{\frac{1}{n}}c_{m},$ and investigate how well this approximates $c_{m}$ as the sample rate grows in relation to $m$. In Section~\ref{fullhalfalg} we examine the multilevel algorithm in the special case where the cosine frequency $m$ is less than half of the initial sample rate. We open by examining the behaviour of the first few iterations of the algorithm when  applied to the  constant (the zero frequency cosine). This provides a case study for the development of the convergence framework for the full algorithm for even functions (a sum of cosines); the proof for odd functions (a sum of sines) follows in the same fashion. The main result shows that the multilevel error can be decomposed into two parts; a truncated trigonometric polynomial term and a remainder term which starts off extremely small (a tolerance we call $\epsilon \approx 10^{-34}$) but which grows at a fixed rate with each iteration. We then investigate the behaviour of the magnitude of truncation term and by monitoring the size of the polynomial coefficients we are able to generate theoretical upper bounds that are very close to those observed numerically. We use these bounds to make an assumption on the behaviour of the magnitude of the truncation terms, which we observe to be correct to an accuracy of $10^{-15}.$  While the growth of the remainder term  may be a theoretical concern, it is not observed numerically because its growth is slow in comparison to the speed at which the truncation part decays; for instance the algorithm would require over $ 60$ iterations before the remainder term becomes larger than $10^{-16}$ and, thus, we see that the algorithm is numerically convergent. In Section~\ref{fullalg}  we examine the general convergence of the special case of the  algorithm which begins by sampling at the integers (i.e., the initial sample rate $=1$). This algorithm is applied to a general cosine $c_{m}$ and we monitor  its behaviour at the early levels, where the sample rate is low relative to $m,$ through to the latter stages when the sample rate overtakes the cosine frequency; at which point we can apply the results from Section~\ref{fullhalfalg}. The results of this section, in their general form, are complicated to express and so in order to aid the reader we illustrate the procedure by following the application of the algorithm to $c_{7},$ which we then use as a template for expressing the conclusions in the general setting. The paper closes in Section~\ref{numerics}  with numerical examples to illuminate the theoretical findings.

\section{Previous Research}\label{overview}

The first work in the RBF community which aimed at mitigating the potential oversampling of information was Floater and Iske \cite{floater}, who proposed a multilevel approximation method where an initial stable approximation is formed on a relatively sparse subset of the data and this is   then refined over multiple levels of residual RBF interpolation on progressively denser subsets.  The original implementation uses Wendland's RBFs (finitely smooth and compactly supported) where the size of the support is scaled to reflect the relative density at a given level. In \cite{hales} a multilevel scheme using polyharmonic splines (finitely smooth and globally supported) on uniform grids was presented and constant reduction in error per level was shown. In \cite{iske}  a modified multilevel method was considered, using thin-plate splines for an initial approximation and with subsequent refinements performed using scaled Wendland RBFs. Wendland and coauthors have explored multilevel schemes using scaled Wendland RBFs for solving both approximation problems and partial differential equations on spheres and compact regions in Euclidean space \cite{farrell,legia1,legia2,wendland}. A hurdle in proving convergence results is that, for the Gaussian RBF, by changing scale of the basis function we also change the underlying approximation space. In relation to this,  we highlight the work of Narcowich et al. \cite{narcowich} who analysed a related scheme but required that sequences of approximation spaces were nested. We emphasise that our approximation spaces are not nested.

As far as the authors are aware the extant  theoretical results on the multilevel method (briefly reviewed in the previous paragraph) apply only to basis functions with finite smoothness. In these cases the numerical stability is improved but one has to accept a saturation point on the accuracy. However, recently  multilevel approximation using scaled Gaussians (infinitely smooth and globally supported) has become of interest due to its key role in multilevel sparse kernel interpolation (MuSIK) and its quasi-interpolatory modification (Q-MuSIK), see \cite{georgoulis,usta}. These approaches have achieved successful results in different areas, see \cite{dong,zhao,usta} for details.

 In this paper we contribute further by presenting a first convergence analysis of the multilevel approximation method  using the Gaussian  basis function. The approach we take  differs from the standard formulation in that we replace interpolation with quasi-interpolation. In order to make the analysis tractable we  will investigate the performance of the scheme when approximating univariate real valued functions with period one. The classical approach to this approximation problem is to use Fourier series, but over the past 50 years, many authors have used shifts of a univariate function \cite{golomb,kushpel,pinkus} and this approach has been adapted to the torus \cite{gomes} and the sphere \cite{wahba}.

  \section{Background and Preliminaries} \label{preliminaries}

Following Delvos \cite{delvos}, we let $\mathcal{C}$ denote the space of continuous real-valued function with period one which we equip with the uniform norm $\|f\|_{\infty}=\sup_{x \in \R}|f(x)|.$ Next we let $\mathcal{L}_{2}$
denote the Hilbert space of square integrable periodic functions with inner product
\[
(f,g):=\int_{0}^{1}f(x)g(x)dx.
\]
The exponentials are given by $e_{k}(x)=\exp(2\pi i k x)$ for $k \in \ZZ.$ The finite Fourier transform of $f \in  \mathcal{L}_{2}$ is given by $\widehat{f}_{k}=(f,e_{-k})$ for $k \in \ZZ$ and its inversion is the Fourier series of $f$ given by $\sum_{k=-\infty}^{\infty}\widehat{f}_{k}e_{k},$
which converges to $f$ in the $\mathcal{L}_{2}$-norm $\|\cdot\|_{2}$ induced by the inner product. Next we let $\mathcal{N}$ denote the space of functions $f \in  \mathcal{L}_{2}$ having absolutely convergent Fourier series, i.e., those for which the norm $\|f\|_{1}=\sum_{k=-\infty}^{\infty}|\widehat{f}_{k}|$
is finite. We have the inclusions $\mathcal{N} \subset \mathcal{C} \subset \mathcal{L}_{2}$ and so for any $f \in \mathcal{N}$ we have the estimates $\|f\|_{2}\le \|f\|_{\infty} \le \|f\|_{1}.$

Next we give some useful results connected to the quasi-interpolant (\ref{quasi}).
First we recall that  the Fourier transform of $\psi$ is $\widehat{\psi} (t)=\exp(-2\pi^2 t^2).$
For a member of the  family of exponentials $e_{m}(x)$ ($m\in \ZZ$) the associated quasi-interpolant can be developed as follows
\[
\begin{aligned}
 Q_{\frac{1}{n}}e_{m}(x):&=
 \sum_{\ell \in \ZZ}e_{m}\left(\frac{\ell}{n}\right)\psi(nx-\ell)=\sum_{j =0}^{n-1}\sum_{\ell \in \ZZ}e_{m}\left(\frac{n\ell+j}{n}\right)\psi(nx-(n\ell+j))
 \\
 &=\sum_{j =0}^{n-1}e_{m}\left(\frac{j}{n}\right)\sum_{\ell \in \ZZ}\psi(n(x-\ell)-j).
 \end{aligned}
\]
Let $\sigma_{j}(x)$ denote the infinite sum appearing in the final line above. We note  that
 $\sigma_{j}(x)$ is $1-$periodic and so we can consider its Fourier expansion
\[
\sigma_{j}(x) = \sum_{\ell \in \ZZ}\psi(n(x-\ell)-j)=\sum_{k=-\infty}^{\infty}\widehat{\sigma}_{k}^{(j)}e_{k}(x)\,\,\,\,{\rm{where}}\,\,\,\,
\widehat{\sigma}_{k}^{(j)}=\int_{0}^{1}\sigma_{j}(x)e_{-k}(x)dx.
\]
Using the periodicity of $\sigma_{j}$ together with an appropriate shift and scale in the variable of integration one can show that the Fourier coefficients are given by:
\[
\widehat{\sigma}_{k}^{(j)} ={1 \over n} e_{-k}\left(\frac{j}{n}\right) \widehat{\psi}\left(\frac{k}{n}\right).
\]
Substituting this back into the expression for $Q_{\frac{1}{n}}e_{m}(x)$ we see that for $m\in \ZZ,$
\EQ{gnexp}{
\begin{aligned}
Q_{\frac{1}{n}} e_m (x) & =  \sum_{j=0}^{n-1} e_{m}\left(\frac{j}{n}\right)\left({1 \over n}\sum_{k=-\infty}^{\infty} e_{-k}\left(\frac{j}{n}\right) \widehat{\psi}\left(\frac{k}{n}\right)e_{k}(x)\right)
\\
& =  \sum_{k=-\infty}^{\infty}\widehat{\psi}\left(\frac{k}{n}\right) e_k(x) \left ( {1 \over n} \sum_{j=0}^{n-1}e_{j}\left(\frac{m-k}{n}\right) \right)
   \\
& =  \sum_{k=-\infty}^{\infty} \widehat{\psi}\left(\frac{nk+m}{n}\right) e_{m+nk} (x)= \sum_{k=-\infty}^{\infty}
\widehat{\psi}\left(k+\frac{m}{n}\right) e_{m+nk} (x).
\end{aligned}}

Our aim is to investigate the convergence rate of our proposed multilevel quasi-interpolation method. To set up the basic framework we shall assume that the target function $f\in \mathcal{N}$ and so possesses a Fourier series $\sum_{k=-\infty}^{\infty}\widehat{f}_{k}e_{k}.$  Then,
using (\ref{gnexp}), the quasi-interpolant is given by
\EQ{gngen}{
Q_{\frac{1}{n}} f  = \sum_{k \in \ZZ} \widehat{f}_{k}Q_{\frac{1}{n}}e_k (x)= \sum_{k=-\infty}^{\infty}
\widehat{f}_{k}\sum_{\ell \in \ZZ} \widehat{\psi}\left(\ell+\frac{k}{n}\right) e_{k+n\ell} (x).}


 We now begin our  investigation by bounding the magnitude of the quasi-interpolant (\ref{gngen}), for which we have
\EQ{initialbound}{
\|Q_{\frac{1}{n}} f\|_{\infty}\le \sum_{k \in \ZZ} \left|\widehat{f}_{k} \right|\sum_{\ell \in \ZZ} \widehat{\psi}\left(\ell+\frac{k}{n}\right)=\sum_{\ell \in \ZZ} \widehat{\psi}\left(\ell+\frac{k}{n}\right)\|f\|_{1}.
}
Following Baxter \cite{baxter}, an application of the Poisson summation formula yields
\[
 \sum_{\ell \in \ZZ} \widehat{\psi}\left(\ell+\frac{k}{n}\right)=\sum_{\ell \in \ZZ} \exp\left(-2\pi^{2}\left(\ell+\frac{k}{n}\right)^{2}\right)=\frac{1}{\sqrt{2\pi}}\sum_{\ell \in \ZZ}e^{-\frac{\ell^{2}}{2}}e^{ \frac{2\ell i\pi k}{n}},
\]
and we observe that this is a theta function of Jacobi type
\EQ{thetadef}{
\theta_{3}(z,q)=\sum_{\ell \in \ZZ}q^{\ell^{2}}e^{2 \ell iz}\,\,\,\,\,q \in \C \,\,\,{\rm{and}}\,\,\, |q|<1.
}
The following product function representation is found in \cite{gradshteyn} (8.181.2)
\EQ{prodreptheta}{
\theta_{3}(z,q)=\prod_{\ell=1}^{\infty}(1+2q^{2\ell-1}\cos(2z)+q^{2(2\ell-1)})(1-q^{2\ell}).
}
If we choose $q = e^{-\frac{1}{2}}$ we can write
\[
\begin{aligned}
E(t):&=\sum_{\ell \in \ZZ} \widehat{\psi}\left(\ell+t\right)=
\frac{1}{\sqrt{2\pi}}\theta_{3}\left(\pi t,e^{-\frac{1}{2}}\right)\\
&=\frac{1}{\sqrt{2\pi}}\prod_{\ell=1}^{\infty}(1+2e^{-\ell+\frac{1}{2}}\cos(2\pi t)+e^{-(2\ell-1)})(1-e^{-\ell}).
\end{aligned}
\]

We observe that $E$ is $1-$periodic and, due to the product representation, it is decreasing on $[0,\frac{1}{2}]$ and increasing on $[\frac{1}{2},1]$ consequently $E$ attains its global max at zero. In view of these observations we can revisit (\ref{initialbound}) and deduce that $\|Q_{\frac{1}{n}} f\|\le \|f\|E(0),$ where, appealing to (\ref{thetadef}), we have
\[
\begin{aligned}
E(0):=\sum_{\ell \in \ZZ}\exp\left(-2\pi^{2}\ell^{2}\right)=\theta_{3}(0,e^{-2\pi^{2}})=1+2e^{-2\pi^{2}}+2e^{-8\pi^{2}}+2e^{-18\pi^{2}}+\ldots
\end{aligned}
\]
where the right hand side are the leading terms in the expansion of (\ref{prodreptheta}) see \cite[16.38.5]{AS}. We can summarise the development above in the following theorem.
 \theo{proposition}{itbound}{
Suppose $f \in \mathcal{N}$. Then, for $n = 1,2,\ldots$,
\[
\|Q_{\frac{1}{n}} f\|_{\infty}\le a\|f\|_{1}\,\,\,\, {\rm{where}} \,\, a = 1+3e^{-2\pi^{2}}=1+3\widehat{\psi}(1).
\]
 Consequently, setting $A = 1+a=2+3\widehat{\psi}(1)$ we have that
$$
\| E_{\frac{1}{n}}f \|_{\infty} = \| Q_{\frac{1}{n}} f-f \|_{\infty} \le A \| f \|_{1}\quad{\rm{and}}\quad \| M_{\frac{1}{n},p} f \| \le A^p \| f \|_{1}.
$$
}

In the radial basis function literature it is usual to develop error estimates for functions belonging to the the appropriate {\sl native space}.  Since we are approximating with Gaussians the correct native space  is
\[
\mathcal{W}_\psi = \left \{ f\in \mathcal{L}_{2}: \| f \|_\psi = \left ( \sum_{k \in \ZZ}^\infty { |\widehat{f}_{k}|^2 \over \widehat \psi(k) } \right )^{1/2} < \infty \right \}.
\]
We observe that, due to the rapid decay of $\widehat \psi(t),$ this space is extremely limited in size and, because of the severity imposed by the native space norm we will, instead, only provide estimates in the wider space $\mathcal{N}$ as defined above.\\

The result of Proposition \ref{itbound} provides a useful general error bound however, in what follows, we will focus attention only on even functions. Thus, as part of this strategy, we will carefully investigate the application of the multilevel algorithm to the cosine family $c_{m}(x)=\cos(2\pi m x),$  $m=\in \Z$ (the same techniques can be used  for the sine family too). To begin the investigation we have the following result which provides two useful identities for the quasi-interpolation of the cosine family.

\theo{Lemma}{quasitrig}{ For $n = 1,2,3,\ldots, $  we have
\[
\begin{aligned}
&(i)\quad Q_{\frac{1}{n}}c_m  =  \sum_{k=-\infty}^\infty \widehat{\psi}\left(k+\frac{m}{n}\right)c_{m+nk}\quad m=0,1,\ldots\\
&(ii)\quad Q_{\frac{1}{n}}c_m=Q_{\frac{1}{n}}c_{m+jn}, \quad j \in \ZZ.
\end{aligned}\]
}
\begin{proof}
 For the first equation, we can use the identity $c_{m}=\frac{1}{2}(e_{m}+e_{-m})$ together with the fact that $\widehat{\psi}(-t)=\widehat{\psi}(t)$ to deduce
\begin{eqnarray*}
Q_{\frac{1}{n}}c_m  & = & {1 \over 2}\left(Q_{\frac{1}{n}}e_m+ Q_{\frac{1}{n}}e_{-m}\right)\\
& = &{1 \over 2} \left(\sum_{k=-\infty}^{\infty}
\widehat{\psi}\left(k+\frac{m}{n}\right) e_{m+nk}+\sum_{k=-\infty}^{\infty}
\widehat{\psi}\left(k-\frac{m}{n}\right) e_{-m+nk}\right)\\
& = &{1 \over 2} \left(\sum_{k=-\infty}^{\infty}
\widehat{\psi}\left(k+\frac{m}{n}\right) e_{m+nk}+\sum_{k=-\infty}^{\infty}
\widehat{\psi}\left(-k-\frac{m}{n}\right) e_{-m-nk}\right)\\
& = &\sum_{k=-\infty}^{\infty}
\widehat{\psi}\left(k+\frac{m}{n}\right)\frac{e_{m+nk}+e_{-m-nk}}{2}=\sum_{k=-\infty}^{\infty}
\widehat{\psi}\left(k+\frac{m}{n}\right)c_{m+nk}.
\end{eqnarray*}
The second equation is an immediate consequence of the first. \end{proof}

We observe  that the second equation of lemma \ref{quasitrig} arises from the aliasing phenomenon i.e., when the sample rate is $n$ then the quasi-interpolant of a cosine of frequency $m$ is indistinguishable from a cosine of frequency $m+jn$ ($j \in \Z).$ It is instructive to consider a sample rate of $n=2^{\ell}$ and write the quasi-interpolant as follows
\EQ{qintcentered}{
Q_{\frac{1}{2^{\ell}}}c_{m}=\widehat{\psi}\left(\frac{m}{2^{\ell}}\right)c_{m}+\sum_{k=1}^{\infty}
\Bigl[\widehat{\psi}\left(\frac{m}{2^{\ell}}-k\right)c_{2^{\ell}k-m}+\widehat{\psi}\left(k+\frac{m}{2^{\ell}}\right)c_{2^{\ell}k+m}\Bigr].
}
Let us consider two cases.\\

$\diamond$ \textbf{Cosine frequency  $\le$   sample rate:} The first case of particular interest is the zero frequency constant function $c_{0}$ where the above formula yields
\[
Q_{\frac{1}{2^{\ell}}}c_{0}=c_{0}+2\widehat{\psi}(1)c_{2^{\ell}}+2\sum_{k=2}^{\infty}\widehat{\psi}(k)c_{2^{\ell}k}.
\]

Due to the very rapid decay of  $\widehat{\psi}(t)$ we can truncate the above series with  hardly any loss in accuracy. Indeed, if we set our tolerance to be $\epsilon =3\widehat{\psi}(2)\sim 10^{-34}$ then we can write
\EQ{Qcons}{
Q_{\frac{1}{2^{\ell}}}c_{0}=c_{0}+2\widehat{\psi}(1)c_{2^{\ell}}+g,\quad {\rm{where}}\,\,\,\,g = 2\sum_{k=2}^{\infty}\widehat{\psi}(k)c_{2^{\ell}k}.
}
We note that the size of remainder term $g$ satisfies
\[
\|g\|_{\infty}\le \|g\|_{1}=2\widehat{\psi}(2)+2\sum_{k=3}^{\infty}\widehat{\psi}(k)\le 3\widehat{\psi}(2)=\epsilon.
\]

Applying the same argument to all other cosine frequencies in the range of consideration we have

\EQ{Qothers}{
Q_{\frac{1}{2^{\ell}}}c_{m}
=\widehat{\psi}\left(\frac{m}{2^{\ell}}-2\right)c_{2^{\ell+1}-m}+\widehat{\psi}\left(\frac{m}{2^{\ell}}-1\right)c_{2^{\ell}-m}+
\widehat{\psi}\left(\frac{m}{2^{\ell}}\right)c_{m}+\widehat{\psi}\left(\frac{m}{2^{\ell}}+1\right)c_{2^{\ell}+m}+g
}
where
\EQ{rem}{\\
g=\widehat{\psi}\left(\frac{m}{2^{\ell}}+2\right)+\sum_{k=3}^{\infty}
\Bigl[\widehat{\psi}\left(\frac{m}{2^{\ell}}-k\right)c_{2^{\ell}k-m}+\widehat{\psi}\left(\frac{m}{2^{\ell}}+k\right)c_{2^{\ell}k+m}\Bigr]
}
denotes the remainder whose magnitude satisfies
\EQ{gbnd}{\\
\|g\|_{\infty}\le\widehat{\psi}\left(2\right)+\widehat{\psi}\left(2\right)+2\sum_{k=3}^{\infty}
\widehat{\psi}\left(k\right)\le 3\widehat{\psi}\left(2\right)=\epsilon.
}


As a side remark we note that for the frequencies in the upper half of the range, i.e., between $2^{\ell-1}$ and $2^{\ell},$  we can  appeal to the aliasing result of Lemma \ref{quasitrig}, which shows that
\EQ{mirror}{
Q_{\frac{1}{2^{\ell}}}c_{2^{\ell-1}+m}=Q_{\frac{1}{2^{\ell}}}c_{2^{\ell-1}-m}, \quad {\rm{for}}\,\,\, m=1,\ldots, 2^{\ell-1}.
}


Using (\ref{Qothers}) we can express the quasi-interpolation error  by
\[
Q_{\frac{1}{2^{\ell}}}c_{m}-c_{m}=\widehat{\psi}\left(\frac{m}{2^{\ell}}-2\right)c_{2^{\ell+1}-m}+\widehat{\psi}\left(\frac{m}{2^{\ell}}-1\right)c_{2^{\ell}-m}+
\left(\widehat{\psi}\left(\frac{m}{2^{\ell}}\right)-1\right)c_{m}+\widehat{\psi}\left(\frac{m}{2^{\ell}}+1\right)c_{2^{\ell}+m}+g.
\]
Furthermore, if we define
\[
F(x):=1+\widehat{\psi}\left(x-2\right)+\widehat{\psi}\left(x-1\right)
-\widehat{\psi}\left(x\right)+\widehat{\psi}\left(x+1\right), \quad x \in [0,1].
\]
then we have the following bound on the quasi-interpolation error
\[
\|Q_{\frac{1}{2^{\ell}}}c_{m}-c_{m}\|_{\infty}\le F\left(\frac{m}{2^{\ell}}\right)+\|g\|_{\infty}\le F\left(\frac{m}{2^{\ell}}\right)+\epsilon.
\]
It is straightforward to verify that the function $F(x)$ is strictly increasing on $[0,1],$ and so, as is to be expected, the approximation error for a particular frequency decreases as the sample rate increases.  The results in Table 1 show the case when the sample rate $=8$ and compares the quasi-interpolation errors of low frequency cosines to their upper bounds based on sampling $F$. Here we see that the actual errors are reasonably close to the upper bounds, especially for the smaller frequencies.

\begin{table}[h]
\caption{Comparison of actual quasi-interpolation errors  with a sample rate $=8$ against the upper bound $F\left(\frac{m}{8}\right)$ for the constant function}
\centering
\begin{tabular}{|c|c|c|} \hline
Frequency  $m$  & $\|Q_{\frac{1}{8}}c_{m}-c_{m}\|_{\infty}$ & $ F\left(\frac{m}{8}\right)$   \\ \hline \hline
0 & 5.35058 (-9) & 5.35058 (-9) \\ \hline
1 & 0.2653969  & 0.2653973   \\ \hline
2 & 0.70877& 0.70880   \\ \hline
3 & 0.93789 & 0.93815  \\ \hline
4 & 0.98561 & 1.00000 \\ \hline
5 & 1.05003  & 1.06184\\ \hline
6 & 1.14906 & 1.29119 \\ \hline
7 & 1.66271 & 1.73460  \\ \hline
8 & 2 - 5.35058 (-9) & 2 + 5.1(-35)  \\ \hline
\end{tabular}
\end{table}

$\diamond$ \textbf{Cosine frequency  $\ge$  sample rate:}

For those frequencies that are integer multiples of the sample rate the aliasing result tells us that the corresponding quasi-interpolants coincide with that of $c_{0},$ i.e.,
\EQ{mults}{
Q_{\frac{1}{2^{\ell}}}c_{r2^{\ell}}=Q_{\frac{1}{2^{\ell}}}c_{0}= c_{0}+2\widehat{\psi}(1)c_{2^{\ell}}+ g.
}
The quasi-interpolation error formula is given by
\[
Q_{\frac{1}{2^{\ell}}}c_{r2^{\ell}}-c_{r2^{\ell}}=1-c_{r2^{\ell}}+2\widehat{\psi}(1)c_{2^{\ell}}+ g.
\]

For the other frequencies of the form $m = r2^{\ell}+n$ where $0<n<2^{\ell}$ we have
\EQ{mods}{
Q_{\frac{1}{2^{\ell}}}c_{m}=Q_{\frac{1}{2^{\ell}}}c_{r2^{\ell}+n}=Q_{\frac{1}{2^{\ell}}}c_{n}=Q_{\frac{1}{2^{\ell}}}c_{m\md 2^{\ell}},
}
and the corresponding quasi interpolation error formula is
\[
\begin{aligned}
&Q_{\frac{1}{2^{\ell}}}c_{r2^{\ell}+n}-c_{r2^{\ell}+n}=Q_{\frac{1}{2^{\ell}}}c_{n}-c_{r2^{\ell}+n}\\&=
\widehat{\psi}\left(\frac{n}{2^{\ell}}-2\right)c_{2^{\ell+1}-n}+\widehat{\psi}\left(\frac{n}{2^{\ell}}-1\right)c_{2^{\ell}-n}+
\widehat{\psi}\left(\frac{n}{2^{\ell}}\right)c_{n}+\widehat{\psi}\left(\frac{n}{2^{\ell}}+1\right)c_{2^{\ell}+n}-c_{r2^{\ell}+n}+g,
\end{aligned}
\]
where, by the nature of the truncation, the size of the remainder term $g$ is bounded by the tolerance $\epsilon.$ If, in this case, we define
\[
G(x)=1+\widehat{\psi}\left(x-2\right)+\widehat{\psi}\left(x-1\right)
+\widehat{\psi}\left(x\right)+\widehat{\psi}\left(x+1\right)=F(x)+2\widehat{\psi}\left(x\right), \quad x \in [0,1],
\]
then  a bound on the size of the quasi-interpolation error is given by

\[
\|Q_{\frac{1}{2^{\ell}}}c_{r2^{\ell}+n}-c_{r2^{\ell}+n}\|_{\infty}\le G\left(\frac{n}{2^{\ell}}\right)+\|g\|_{\infty}\le G\left(\frac{n}{2^{\ell}}\right)+\epsilon.
\]

We observe that $G(0)=G(1)=2+2\widehat{\psi}\left(1\right)+\order(\epsilon).$ Further, it can easily be shown that $G$ is decreasing on $[0,\frac{1}{2}],$ increasing on $[\frac{1}{2},1]$ and, at $x=\frac{1}{2},$ it attains its minimum value
\[
G\left( \frac{1}{2}\right)=1+2\widehat{\psi}\left(\frac{1}{2}\right)+2\widehat{\psi}\left(\frac{3}{2}\right).
\]
The results in Table 2 show the case when the sample rate $=8$ and compares the quasi-interpolation errors of higher frequency cosines to their upper bounds, based on sampling $G.$. Here we see, again,  that the actual errors are reasonably close to the upper bounds.
\begin{table}[h]
\caption{Comparison of actual quasi-interpolation errors with sample rate $=8$ of relatively high frequency cosines $c_{8+n}$ against the upper bound $G\left(\frac{n}{8}\right)$ for the constant function}
\centering
\begin{tabular}{|c|c|c|} \hline
Frequency  $8+n$  & $\|Q_{\frac{1}{8}}c_{8+n}-c_{8+n}\|$ & $ G\left(\frac{n}{8}\right)$   \\ \hline \hline
9 & 1.69068 & 1.73460 \\ \hline
10 & 1.23559  & 1.29122  \\ \hline
11 & 1.05881 & 1.06274   \\ \hline
12 & 1.00720 & 1.01438\\ \hline
13 & 1.06037 & 1.06274 \\ \hline
14 & 1.26152 & 1.29122  \\ \hline
15 & 1.71737 & 1.73460 \\ \hline
16 & 1.99998 & 2.0000 \\ \hline
\end{tabular}
\end{table}

\section{Convergence of the discrete algorithm for a low frequency cosine} \label{fullhalfalg}

In this section we will assume that the starting sample rate of the multilevel algorithm is $2^{\ell}$ ($\ell\ge1$) and that it is
 applied to
 $c_{m}$ the cosine function with a fixed frequency $m<2^{\ell-1}.$ In order to shed light on the approach that we will take we begin by developing the first few iterations of the method applied to the constant function, this will act as a prototype for approximation of all frequencies.\\

$\diamond$ \textbf{Approximation of the constant function:} We begin by approximating the constant function via  quasi-interpolation at the integers, this coincides with the first level of the algorithm. In what follows we will frequently appeal to (\ref{Qcons}), for the quasi-interpolation formula, but we will ignore the remainder terms and any other terms that develop in the course of the iterations whose  magnitudes are of order $\epsilon.$ We begin at level one where we have

\beqns
M_{1,1} c_0 & = & Q_{1}c_{0}-c_{0}=2\widehat{\psi}(1)c_{1}+ \order(\epsilon),
\eeqns

In what follows we will ignore all terms that are smaller or of the same order as $\epsilon=3\hpsi(2) \approx 10^{-34}.$ The multilevel error at the next level is given by
\EQ{connect-one}{
M_{1,2}c_{0}=Q_{\frac{1}{2}}M_{1,1} c_0-M_{1,1} c_0\approx 2\widehat{\psi}(1)\left(Q_{\frac{1}{2}}c_{1}-c_{1}\right)
}
and appealing to (\ref{mirror}) we have
\[
M_{1,2}c_{0}\approx 2\widehat{\psi}(1)\left(\left(2\widehat{\psi}\left(\frac{1}{2}\right)-1\right)c_{1}+2\widehat{\psi}\left(\frac{3}{2}\right)c_{3}\right).
\]

Moving to the next level we have
\[
M_{1,3}c_{0}=Q_{\frac{1}{4}}M_{1,2} c_0-M_{1,2} c_0\approx 2\widehat{\psi}(1)\left(\left(2\widehat{\psi}\left(\frac{1}{2}\right)-1\right)\left(Q_{\frac{1}{4}}c_{1}-c_{1}\right)
+2\widehat{\psi}\left(\frac{3}{2}\right)\left(Q_{\frac{1}{4}}c_{3}-c_{3}\right)\right),
\]
and appealing to (\ref{Qothers}) and the aliasing property we have
\[
 Q_{\frac{1}{4}}c_{1}=Q_{\frac{1}{4}}c_{3}\approx
\widehat{\psi}\left(\frac{1}{4}\right)c_{1}+\widehat{\psi}\left(\frac{3}{4}\right)c_{3}+
\widehat{\psi}\left(\frac{5}{4}\right)c_{5}+\widehat{\psi}\left(\frac{7}{4}\right)c_{7}.
\]

Substituting and simplifying we have
\[
\begin{aligned}
M_{1,3}c_{0}\approx 2\widehat{\psi}(1)\Bigl[&\Bigl\{\left(2\widehat{\psi}\left(\frac{1}{2}\right)-1\right)\left(\widehat{\psi}\left(\frac{1}{4}\right)-1\right)
+2\widehat{\psi}\left(\frac{3}{2}\right)\widehat{\psi}\left(\frac{1}{4}\right)\Bigr\}c_{1}\\
+&2\widehat{\psi}\left(\frac{7}{4}\right)\Bigl\{2\widehat{\psi}\left(\frac{1}{2}\right)+2\widehat{\psi}\left(\frac{3}{2}\right)-1\Bigr\}c_{7}\\
+&\Bigl\{\left(2\widehat{\psi}\left(\frac{1}{2}\right)-1\right)\widehat{\psi}\left(\frac{3}{4}\right)
+2\widehat{\psi}\left(\frac{3}{2}\right)\left(\widehat{\psi}\left(\frac{3}{4}\right)-1\right)\Bigr\}c_{3}\\
+&2\widehat{\psi}\left(\frac{5}{4}\right)\Bigl\{2\widehat{\psi}\left(\frac{1}{2}\right)+2\widehat{\psi}\left(\frac{3}{2}\right)-1\Bigr\}c_{5}
\Bigr].
\end{aligned}
\] We have written the frequencies in this order because we wish to pair up the frequencies which are aliased in the next iteration of the algorithm: here $Q_{\frac{1}{8}}c_{1}=Q_{\frac{1}{8}}c_{7}$ and $Q_{\frac{1}{8}}c_{3}=Q_{\frac{1}{8}}c_{5}$.
By ignoring the terms in the above expression that are of order $\epsilon$ we  have
 \[
\begin{aligned}
M_{1,3}c_{0}\approx 2\widehat{\psi}(1)\left(2\widehat{\psi}\left(\frac{1}{2}\right)-1\right)
&\Bigl[\left(\widehat{\psi}\left(\frac{1}{4}\right)-1\right)c_{1}
+\widehat{\psi}\left(\frac{3}{4}\right)c_{3}+2\widehat{\psi}\left(\frac{5}{4}\right)c_{5}\Bigr]\\
+4\widehat{\psi}(1)\widehat{\psi}\left(\frac{3}{2}\right)&\Bigl[\widehat{\psi}\left(\frac{1}{4}\right)c_{1}
+\left(\widehat{\psi}\left(\frac{3}{4}\right)-1\right)c_{3}\Bigr].
\end{aligned}
\]

Let us look at one more iteration of the algorithm, we have
 \[
\begin{aligned}
M_{1,4}c_{0}=Q_{\frac{1}{8}}M_{1,3} c_0-M_{1,3} c_0 &\\
\approx 2\widehat{\psi}(1)\left(2\widehat{\psi}\left(\frac{1}{2}\right)-1\right)&\Bigl[\left(\widehat{\psi}\left(\frac{1}{4}\right)-1\right)(Q_{\frac{1}{8}}c_{1}-c_{1})
+\widehat{\psi}\left(\frac{3}{4}\right)(Q_{\frac{1}{8}}c_{3}-c_{3})+2\widehat{\psi}\left(\frac{5}{4}\right)(Q_{\frac{1}{8}}c_{5}-c_{5})\Bigr]\\
+4\widehat{\psi}(1)\widehat{\psi}\left(\frac{3}{2}\right)&\Bigl[\widehat{\psi}\left(\frac{1}{4}\right)(Q_{\frac{1}{8}}c_{1}-c_{1})
+\left(\widehat{\psi}\left(\frac{3}{4}\right)-1\right)(Q_{\frac{1}{8}}c_{3}-c_{3})\Bigr].
\end{aligned}
\]
Once again, appealing to (\ref{Qothers}) and the aliasing property we have
\beqns
Q_{\frac{1}{8}} c_1 & = & \hpsi\left(\frac{1}{8}\right) c_1 + \hpsi\left(\frac{7}{8}\right) c_7 + \hpsi\left(\frac{9}{8}\right) c_9 + \hpsi\left(\frac{15}{8}\right) c_{15} \\
Q_{\frac{1}{8}} c_3 =Q_{\frac{1}{8}} c_5 &=&\hpsi\left(\frac{3}{8}\right) c_3  + \hpsi\left(\frac{5}{8}\right) c_5 + \hpsi\left(\frac{11}{8}\right) c_{11} + \hpsi\left(\frac{13}{8}\right) c_{13}.
\eeqns

Substituting, simplifying and ignoring  terms that are of order $\epsilon$ we  have
\[
\begin{aligned}
M_{1,4} c_0  \approx\,\,  &2\widehat{\psi}(1)\Bigl[\left(2\widehat{\psi}\left(\frac{1}{2}\right)-1\right)\left(\widehat{\psi}\left(\frac{1}{4}\right)-1\right)
+\hpsi\left(\frac{3}{2}\right)\hpsi\left(\frac{1}{4}\right)\Bigr]\left(\widehat{\psi}\left(\frac{1}{8}\right)-1\right)c_{1}\\
+&2\widehat{\psi}(1)\Bigl[\left(2\widehat{\psi}\left(\frac{1}{2}\right)-1\right)\left(\widehat{\psi}\left(\frac{1}{4}\right)-1\right)\Bigr]
\left(\widehat{\psi}\left(\frac{7}{8}\right)c_{7}+\widehat{\psi}\left(\frac{9}{8}\right)c_{9}\right)\\
+&2\widehat{\psi}(1)\Bigl[\left(2\widehat{\psi}\left(\frac{1}{2}\right)-1\right)\widehat{\psi}\left(\frac{3}{4}\right)+
+2\hpsi\left(\frac{3}{2}\right)\left(\hpsi\left(\frac{3}{4}\right)-1\right)\Bigr]\Bigl\{\left(\widehat{\psi}\left(\frac{3}{8}\right)-1\right)c_{3}+
\widehat{\psi}\left(\frac{5}{8}\right)c_{5}\Bigr\}\\
+&4\widehat{\psi}(1)\left(2\widehat{\psi}\left(\frac{1}{2}\right)-1\right)\hpsi\left(\frac{5}{4}\right)\Bigl\{\widehat{\psi}
\left(\frac{3}{8}\right)c_{3}+
\left(\widehat{\psi}\left(\frac{5}{8}\right)-1\right)c_{5}\Bigr\}\\
+&2\widehat{\psi}(1)\left(2\widehat{\psi}\left(\frac{1}{2}\right)-1\right)\hpsi\left(\frac{5}{4}\right)\hpsi\left(\frac{11}{8}\right)c_{11}
\end{aligned}
\]
again pairing frequencies which alias at the next iteration. We conclude by expressing the error expressions with the approximate numerical values (using the notation $a(-n)=a\times 10^{-n}$) of the coefficients, again ignoring those whose magnitude is of the order $\epsilon.$
\[
\begin{aligned}
M_{1,1}c_{0}&\approx 5.4(-9) c_{1}\\
M_{1,2}c_{0} &\approx -5.3(-9) c_{1}+5.5(-28) c_{3}\\
M_{1,3}c_{0}&\approx 3.7(-9) c_{1}-7.9(-14) c_{3}- 2(-22) c_{5}\\
M_{1,4}c_{0}&\approx -9.9(-10) c_{1}+7.4(-14) c_{3}- 3.5(-17) c_{5}+1(-15) c_{7}+5.3(-20) c_{9}-4.9(-30) c_{11}\\
M_{1,5}c_{0}&\approx 7.3(-11) c_{1}-3.7(-14) c_{3}+3(-17) c_{5}-1(-15) c_{7}+1.9(-18) c_{9}+3.2(-21) c_{11}\\
&+1.6(-19) c_{13}-2.9(-17) c_{15}-2(-19) c_{17}-6(-26) c_{19}-6(-32) c_{21}+1(-32) c_{23}
\end{aligned}
\]
Continuing the development, but displaying only the terms which are relatively numerically significant, we have

\[
\begin{aligned}
M_{1,6}c_{0}&\approx -1.4(-12) c_{1}-5.9(-15) c_{3}+6.1(-16)c_{7}-1.2(-17) c_{5}+2.9(-17) c_{15}-1.5(-18) c_{9}\\
M_{1,7}c_{0}&\approx 6.8(-15) c_{1}-2.5(-16) c_{3}-1.3(-16)c_{7}-1.9(-17) c_{15}+1.3(-18) c_{5}-5(-19) c_{9}
\end{aligned}
\]

Hopefully this example gives the reader an idea of how the iterations proceed. There are two important  points to highlight. The first is to acknowledge that
the number of terms in the approximation grows at each level. For the constant function, when the sample rate is $1$ the first level approximation involves only $c_{1}$.  When the sample rate doubles to $2$ the second level approximation refines the previous $c_{1}$ contribution  and introduces a new contribution from $c_{3}$.  When the sample rate doubles again to $4$ the third level approximation refines the previous contributions (from $c_{1}$ and $c_{3}$) and  introduces new contributions from $c_{5}$ and $c_{7}.$ The algorithm continues in this way at every iteration. Clearly the error expansion grows rapidly but, as we see in the example, we are able to ignore those terms whose contributions are smaller than the tolerance $\epsilon \approx 10^{-34}.$ This leads to the second point which is to note
 that the error is determined by the dominant coefficient(s) in the error expression. By examining the multilevel error expressions for the constant function this suggest that we define the sequence
 \[
 \delta_{1}=2\widehat{\psi}(1), \quad \delta_{2}=2\widehat{\psi}(1)\left(2\widehat{\psi}\left(\frac{1}{2}\right)-1\right),\,\,
 \delta_{p}=\delta_{2}\prod_{j=3}^{p}\left(\widehat{\psi}\left(\frac{1}{2^{j-1}}\right)-1\right)
 \,\, p\ge3.
 \]
 The sequence $\delta_{p}$ enters the error expressions  (at level $p$) as coefficients of $c_{1}.$ A comparison of
 $\delta_{p}$ against $\|M_{1,p} c_0\|_{\infty}$ the realised multilevel error at the $p^{th}$ iteration  is displayed in the table below. Here we observe that initially the sequence tracks the true error extremely well. However, we notice that after the $6^{th}$ iteration the coefficient $\delta_{p}$ is no longer dominant contribution in the error expression.


\begin{table}[h]
\caption{Comparison of level $p$ numerical errors against $\delta_{p}$ for the constant function}
\centering
\begin{tabular}{|c|c|c|} \hline
Level $p$  & $\| M_{1,p} c_0\|_\infty$ & $|\delta_{p}|$   \\ \hline \hline
1 & 5.35058 (-9) & 5.35058 (-9) \\ \hline
2 & 5.27361 (-9)  & 5.27361 (-9)    \\ \hline
3 & 3.73779 (-9)  & 3.73787 (-9)    \\ \hline
4 & 9.91944 (-10) & 9.92020 (-10)   \\ \hline
5 & 7.35780 (-11) & 7.36163 (-11)  \\ \hline
6 & 1.39920 (-12) & 1.40548 (-12) \\ \hline
7 & 6.48407 (-15) & 8.13574 (-18)   \\ \hline
8 & 4.61947 (-16) & 2.45009 (-21)   \\ \hline
9 & 4.20474 (-16) & 1.84482 (-25)   \\ \hline
\end{tabular}
\end{table}

We now devote attention to a more general framework where multilevel error bounds are derived for $c_{m}$ the cosine function of frequency $m$ which we assume to be less than half of the sample rate, i.e., $m < 2^{\ell-1}.$  We will continue to respect $\epsilon = 3\hpsi\left(2\right)\approx 10^{-34}$ as the minimum tolerance and will develop a formula for the multilevel error that aggregates the $\order(\epsilon)$ factors separately from the other contributions.\\
%
%
%



 In the example for the constant function we saw how  the expressions for the cosine coefficients appearing in the error expansions become increasingly more complex as the multilevel algorithm progresses. In order to circumvent this we will develop a recursive relation for the  coefficients, so that their values at a given level are expressed neatly in terms of those from the previous level.\\

$\diamond$ \textbf{Level one:} Our aim throughout this part of the analysis, as with the example of approximating the constant function, is to express the multilevel error in two parts, a main truncation  consisting of the terms that contribute more than the tolerance $\epsilon$ and a remainder term whose size is less than $\epsilon.$ For the first level we write
\EQ{expone}{
M_{\frac{1}{2^{\ell}},1}c_{m}=Q_{\frac{1}{2^{\ell}}}c_{m}-c_{m}=T_{\frac{1}{2^{\ell}},1}c_{m} +g^{(1)},
}
where $T_{\frac{1}{2^{\ell}},1}c_{m}$ is the level one truncation which, using (\ref{Qothers}), is given by
\[
T_{\frac{1}{2^{\ell}},1}c_{m}=\overline{\alpha}_{0}^{(1)}c_{2^{\ell+1}-m}+ \alpha_{0}^{(1)}c_{m}+\overline{\alpha}_{1}^{(1)}c_{2^{\ell}-m}+\alpha_{1}^{(1)}c_{2^{\ell}+m},
\]
whose coefficients are
\EQ{startcoefs}{
\overline{\alpha}_{0}^{(1)}=\widehat{\psi}\left(\frac{m}{2^{\ell}}-2\right),\,\,\quad \alpha_{0}^{(1)}=\widehat{\psi}\left(\frac{m}{2^{\ell}}\right)-1,\quad\overline{\alpha}_{1}^{(1)}=\widehat{\psi}\left(\frac{m}{2^{\ell}}-1\right),
\,\,\quad \alpha_{1}^{(1)}=\widehat{\psi}\left(\frac{m}{2^{\ell}}+1\right).
}

The remainder term $g^{(1)}$ coincides precisely with $g$ in  (\ref{rem}) and thus, appealing to (\ref{gbnd}), we have that $\|g^{(1)}\|_{\infty}\le \epsilon.$  

We remark that the notation used to represent the coefficients in the truncation expression has been deliberately chosen to highlight the pairs of cosines whose quasi-interpolants at next level will match due to aliasing. In particular, for this case we have
\EQ{ale}{
\begin{aligned}
{\rm{for}}\,\,{\rm{sub-indices}}\,&=\,0: \,\,\,Q_{\frac{1}{2^{\ell+1}}}c_{2^{\ell+1}-m}=Q_{\frac{1}{2^{\ell+1}}}c_{m};\\
{\rm{for}}\,\,{\rm{sub-indices}}\,&=\,1: \,\,\,Q_{\frac{1}{2^{\ell+1}}}c_{2^{\ell}-m}=Q_{\frac{1}{2^{\ell}+1}}c_{2^{\ell}+m}.
\end{aligned}
}

$\diamond$ \textbf{Level two}. At the second level the number of sample points is doubled to $2^{\ell+1}$ and multilevel error is given by
\[
\begin{aligned}
M_{\frac{1}{2^{\ell}},2}c_{m}&=Q_{\frac{1}{2^{\ell+1}}}M_{\frac{1}{2^{\ell}},1}c_{m}-M_{\frac{1}{2^{\ell}},1}c_{m}
=Q_{\frac{1}{2^{\ell+1}}}T_{\frac{1}{2^{\ell}},1}c_{m}-T_{\frac{1}{2^{\ell}},1}c_{m}+Q_{\frac{1}{2^{\ell+1}}}g_{1}-g_{1}.
\end{aligned}
\]

Focussing on the truncation term we have
\[
\begin{aligned}
Q_{\frac{1}{2^{\ell+1}}}T_{\frac{1}{2^{\ell}},1}c_{m}-T_{\frac{1}{2^{\ell}},1}c_{m}
&=\overline{\alpha}_{0}^{(1)}\left(Q_{\frac{1}{2^{\ell+1}}}c_{2^{\ell+1}-m}-c_{2^{\ell+1}-m}\right)+ \alpha_{0}^{(1)}\left(Q_{\frac{1}{2^{\ell+1}}}c_{m}-c_{m}\right)\\
&+\overline{\alpha}_{1}^{(1)}\left(Q_{\frac{1}{2^{\ell+1}}}c_{2^{\ell}-m}-c_{2^{\ell}-m}\right)
+\alpha_{1}^{(1)}\left(Q_{\frac{1}{2^{\ell+1}}}c_{2^{\ell}+m}-c_{2^{\ell}+m}\right)\\
&=\overline{\alpha}_{0}^{(1)}\left(Q_{\frac{1}{2^{\ell+1}}}c_{m}-c_{2^{\ell+1}-m}\right)+ \alpha_{0}^{(1)}\left(Q_{\frac{1}{2^{\ell+1}}}c_{m}-c_{m}\right)\\
&+\overline{\alpha}_{1}^{(1)}\left(Q_{\frac{1}{2^{\ell+1}}}c_{2^{\ell}+m}-c_{2^{\ell}-m}\right)
+\alpha_{1}^{(1)}\left(Q_{\frac{1}{2^{\ell+1}}}c_{2^{\ell}+m}-c_{2^{\ell}+m}\right),
\end{aligned}
\]
where we have used the aliasing formulae (\ref{ale}). Now computing the errors for the quasi-interpolants we have, for the first pair, that
\[
\begin{aligned}
Q_{\frac{1}{2^{\ell+1}}}c_{m}-c_{m}&=\widehat{\psi}\left(\frac{m}{2^{\ell+1}}-2\right)c_{2^{\ell+2}-m}
+\widehat{\psi}\left(\frac{m}{2^{\ell+1}}-1\right)c_{2^{\ell+1}-m}\\
&+
\left(\widehat{\psi}\left(\frac{m}{2^{\ell+1}}\right)-1\right)c_{m}+\widehat{\psi}\left(\frac{m}{2^{\ell+1}}+1\right)c_{2^{\ell+1}+m}+g_{0}^{(2)}
\end{aligned}
\]
and
\[
\begin{aligned}
Q_{\frac{1}{2^{\ell+1}}}c_{m}-c_{2^{\ell+1}-m}&=\widehat{\psi}\left(\frac{m}{2^{\ell+1}}-2\right)c_{2^{\ell+2}-m}
+\left(\widehat{\psi}\left(\frac{m}{2^{\ell+1}}-1\right)-1\right)c_{2^{\ell+1}-m}\\&+
\widehat{\psi}\left(\frac{m}{2^{\ell+1}}\right)c_{m}+\widehat{\psi}\left(\frac{m}{2^{\ell+1}}-1\right)c_{2^{\ell+1}+m}+g_{0}^{(2)}
\end{aligned}
\]
where, in both cases, we can appeal to the previous section, specifically (\ref{Qothers})-(\ref{gbnd}), to deduce that the remainder term
$g_{0}^{(2)}$ satisfies $\|g_{0}^{(2)}\|_{\infty}\le \epsilon.$
For the second pair we have
\[
\begin{aligned}
Q_{\frac{1}{2^{\ell+1}}}c_{2^{\ell}-m}-c_{2^{\ell}-m}&=\widehat{\psi}\left(\frac{m}{2^{\ell+1}}-\frac{3}{2}\right)c_{3\cdot2^{\ell}-m}
+\left(\widehat{\psi}\left(\frac{m}{2^{\ell+1}}-\frac{1}{2}\right)-1\right)c_{2^{\ell}-m}\\&+
\widehat{\psi}\left(\frac{m}{2^{\ell+1}}+\frac{1}{2}\right)c_{2^{\ell}+m}+
\widehat{\psi}\left(\frac{m}{2^{\ell+1}}+\frac{3}{2}\right)c_{3\cdot2^{\ell}+m}+g_{1}^{(2)}
\end{aligned}
\]
and
\[
\begin{aligned}
Q_{\frac{1}{2^{\ell+1}}}c_{2^{\ell}+m}-c_{2^{\ell}+m}&=\widehat{\psi}\left(\frac{m}{2^{\ell+1}}-\frac{3}{2}\right)c_{3\cdot2^{\ell}-m}
+\widehat{\psi}\left(\frac{m}{2^{\ell+1}}-\frac{1}{2}\right)c_{2^{\ell}-m}\\&+
\left(\widehat{\psi}\left(\frac{m}{2^{\ell+1}}+\frac{1}{2}\right)-1\right)c_{2^{\ell}+m}+
\widehat{\psi}\left(\frac{m}{2^{\ell+1}}+\frac{3}{2}\right)c_{3\cdot2^{\ell}+m}+g_{1}^{(2)}
\end{aligned}
\]
where, as above, the remainder term
$g_{1}^{(2)}$ satisfies $\|g_{1}^{(2)}\|_{\infty}\le \epsilon.$
%
With this we can develop the overall level two error expression as
\[
M_{\frac{1}{2^{\ell}},2}c_{m}=T_{\frac{1}{2^{\ell}},2}c_{m}+g^{(2)},
\]
where the remainder term is given by
\EQ{gtwo}{
g^{(2)}=(\overline{\alpha}_{0}^{(1)}+\alpha_{0}^{(1)})g_{0}^{(2)}+(\overline{\alpha}_{1}^{(1)}+\alpha_{1}^{(1)})g_{1}^{(2)}+Q_{\frac{1}{2^{\ell+1}}}g_{1}-g_{1}.
}
and where the truncation term has the form
\EQ{Ttwo}
{
\begin{aligned}
T_{\frac{1}{2^{\ell}},2}c_{m}&=
\overline{\alpha}_{0}^{(2)}c_{2^{\ell+2}-m}&+&\overline{\alpha}_{1}^{(2)}c_{3\cdot2^{\ell}-m}
&+&\overline{\alpha}_{2}^{(2)}c_{\cdot2^{\ell+1}-m}&+&\overline{\alpha}_{3}^{(2)}c_{\cdot2^{\ell}-m}\\&
+\alpha_{0}^{(2)}c_{m}
&+&\alpha_{1}^{(2)}c_{2^{\ell}+m}
&+&\alpha_{2}^{(2)}c_{2^{\ell+1}+m}
&+&\alpha_{3}^{(2)}c_{3\cdot2^{\ell}+m}.
\end{aligned}
}
Inspecting the error expansion, the coefficients introduced above are
\[
\begin{aligned}
\overline{\alpha}_{j}^{(2)}&=(\overline{\alpha}_{j}^{(1)}+\alpha_{j}^{(1)})
\widehat{\psi}\left(2-\frac{m}{2^{\ell+1}}-\frac{j}{2}\right)\,\,\,\, (j=0,1),\\
\overline{\alpha}_{2^{1}+j}^{(2)}&=\overline{\alpha}_{j}^{(1)}\Bigl[\widehat{\psi}\left(1-\frac{m}{2^{\ell+1}}-\frac{j}{2}\right)-1\Bigr]
+
\alpha_{j}^{(1)}\widehat{\psi}\left(1-\frac{m}{2^{\ell+1}}-\frac{j}{2}\right)\,\,\,\, (j=0,1),\\
\alpha_{j}^{(2)}&=\alpha_{j}^{(1)}\Bigl[\widehat{\psi}\left(\frac{m}{2^{\ell+1}}+\frac{j}{2}\right)-1\Bigr]+
\overline{\alpha}_{j}^{(1)}\widehat{\psi}\left(\frac{m}{2^{\ell+1}}+\frac{j}{2}\right)\,\,\,\, (j=0,1),\\
\alpha_{2^{1}+j}^{(2)}&=(\overline{\alpha}_{j}^{(1)}+\alpha_{j}^{(1)})\widehat{\psi}\left(1+\frac{m}{2^{\ell+1}}+\frac{j}{2}\right)\,\,\,\, (j=0,1).\\
\end{aligned}
\]

To investigate the size of the remainder at level two we can proceed as follows.
\[
\begin{aligned}
\|g^{(2)}\|_{\infty}&\le \|g^{(1)}-Q_{\frac{1}{2^{\ell+1}}}g^{(1)}\|_{\infty}+\left(|\overline{\alpha}_{0}^{(1)}|+|\alpha_{0}^{(1)}|\right)\|g_{0}^{(1)}\|
+\left(|\overline{\alpha}_{1}^{(1)}|+|\alpha_{1}^{(1)}|\right)\|g_{1}^{(1)}\|_{\infty}\\
&\le \|E_{\frac{1}{2^{\ell+1}}}g^{(1)}\|_{\infty}+\left(|\overline{\alpha}_{0}^{(1)}|+|\alpha_{0}^{(1)}|
+|\overline{\alpha}_{1}^{(1)}|+|\alpha_{1}^{(1)}|\right)\epsilon.
\end{aligned}
\]

We can evoke Proposition \ref{itbound}  for the following crude bound on the quasi-interpolation error of $g_{1}$
\[
\|E_{\frac{1}{2^{\ell+1}}}g^{(1)}\|_{\infty}\le A\|g^{(1)}\|_{\infty}\le A\epsilon=(2+ 3\widehat{\psi}(1))\epsilon.
\]
In addition, using (\ref{startcoefs}), we also have that
\[
|\overline{\alpha}_{0}^{(1)}|+|\alpha_{0}^{(1)}|
+|\overline{\alpha}_{1}^{(1)}|+|\alpha_{1}^{(1)}|=\widehat{\psi}\left(\frac{m}{2^{\ell}}-2\right)+1-\widehat{\psi}\left(\frac{m}{2^{\ell}}\right)
+\widehat{\psi}\left(\frac{m}{2^{\ell}}-1\right)+\widehat{\psi}\left(\frac{m}{2^{\ell}}+1\right)\le 1+3\widehat{\psi}\left(\frac{1}{2}\right).
\]
Taking these bounds into consideration we can conclude that
\[
\|g^{(2)}\|_{\infty}\le 2A\epsilon.
\]


The preceding analysis of the first two levels provides sufficient insight to establish the following, more general result.
\begin{proposition} \label{propiter}
Let $m$ be the fixed frequency of the cosine $c_{m}$ and assume $m < 2^{\ell-1}$, for some $\ell \ge 2$. Then
\beqns
M_{\frac{1}{2^{\ell}},p} c_m & =  & T_{\frac{1}{2^{\ell}},p} c_m + g_p,\quad {\rm{where}}\quad
 T_{\frac{1}{2^{\ell}},p} c_m  =  \sum_{j=0}^{2^p-1} ( \overline{\alpha}_{j}^{(p)}  c_{m-2^{\ell}(2^{p}-j)} +\alpha_{j}^{(p)} c_{m + 2^{\ell}j}),
\eeqns
where $\| g_p \|_{\infty} \le pA^{p-1} \epsilon$. The truncation coefficients are defined recursively.
Specifically, for $p=1$ the initial coefficients are given by (\ref{startcoefs}), and then for $p>1$ we have

 \EQ{coeffs}{
\begin{aligned}
\overline{\alpha}_{j}^{(p+1)}&=(\overline{\alpha}_{j}^{(p)}+\alpha_{j}^{(p)})
\widehat{\psi}\left(2-\frac{m}{2^{\ell+p}}-\frac{j}{2^{p}}\right),\,\,\,\, &j&=0,1,\ldots,2^{p}-1;\\
\overline{\alpha}_{2^{p}+j}^{(p+1)}&=\overline{\alpha}_{j}^{(p)}
\Bigl[\widehat{\psi}\left(1-\frac{m}{2^{\ell+p}}-\frac{j}{2^{p}}\right)-1\Bigr]
+
\alpha_{j}^{(p)}\widehat{\psi}\left(1-\frac{m}{2^{\ell+p}}-\frac{j}{2^{p}}\right),\,\,\,\, &j&=0,1,\ldots,2^{p}-1;\\
\alpha_{j}^{(p+1)}&=\alpha_{j}^{(p)}\Bigl[\widehat{\psi}\left(\frac{m}{2^{\ell+p}}+\frac{j}{2^{p}}\right)-1\Bigr]+
\overline{\alpha}_{j}^{(p)}\widehat{\psi}\left(\frac{m}{2^{\ell+p}}+\frac{j}{2^{p}}\right),\,\,\,\, &j&=0,1,\ldots,2^{p}-1;\\
\alpha_{2^{p}+j}^{(p+1)}&=(\overline{\alpha}_{j}^{(p)}+\alpha_{j}^{(p)})\widehat{\psi}\left(1+\frac{m}{2^{\ell+p}}+\frac{j}{2^{p}}\right),\,\,\,\, &j&=0,1,\ldots,2^{p}-1.
\end{aligned}
}

\end{proposition}

\begin{proof} The result can be established by induction on $p.$ Indeed, the preceding error analysis for the level one case establishes the result for $p=1.$ Assuming the result for a general level $p$ (the inductive hypothesis) one can then  mimic the  methodology of the level two analysis to inductively establish the stated result for the $p+1$ level.
 \end{proof}

 We alert the reader to the fact that the constant $A,$ appearing in the bound for the general remainder term $g^{(p)},$ is  greater than $2.$ This is significant because it clearly rules out the possibility of concluding that the multilevel algorithm will converge as $p\to \infty.$ However, from a numerical perspective, the success of the algorithm will be judged on how quickly multilevel error achieves a desired level of accuracy.  If, for instance, we target an accuracy of $10^{-15}$ for the approximation error, as was observed in the example of approximating the constant,  then this will be achieved if the size of the truncation $\|T_{\frac{1}{2^{\ell}},p} c_m \|_{\infty}$ reaches $\sim 10^{-15}$ before the size of the remainder term grows higher than  $\sim 10^{-15};$ based on the bound we have for $\|g^{(p)}\|_{\infty}$ this would happen at the $60^{th}$ iteration of the algorithm. In the case of approximating the constant
 we saw that  this level of accuracy was achieved after only $7$ iteration and that  the product sequence $(\delta_{p})_{p\ge 0}$ tracked the numerical error extremely well across the early levels of the algorithms but that its influence fades when the numerical error is governed by other cosine frequencies that enter the error expression. \\

 In what follows we will use  Proposition~\ref{propiter} to derive an analytic sequence that tracks the error  when the cosine  frequency is less than half of the sampling rate.  To set the scene we consider the following partition
 \[
 [1,2^{\ell-1}]=[2^{0},2^{1})\cup [2^{1},2^{2})\cup\cdots\cup [2^{\ell-2},2^{\ell-1}).
 \]

 Suppose $m\in [2^{k-1},2^{k})$ then $\frac{1}{2^{1+\ell-k}}\le \frac{m}{2^{\ell}}<\frac{1}{2^{\ell-k}}$ and so, using (\ref{coeffs}), we can bound the coefficients appearing in the error estimate at level $p+1$ as multiples of those from the previous level. Specifically we have that
  \EQ{boundedcoeffs}{
\begin{aligned}
|\overline{\alpha}_{j}^{(p+1)}|&\le(|\overline{\alpha}_{j}^{(p)}|+|\alpha_{j}^{(p)}|)
\widehat{\psi}\left(2-\frac{1}{2^{\ell-k+p}}-\frac{j}{2^{p}}\right),\,\,\,\, &j&=0,1,\ldots,2^{p}-1;\\
|\overline{\alpha}_{2^{p}+j}^{(p+1)}|&\le|\overline{\alpha}_{j}^{(p)}|
\Bigl[1-\widehat{\psi}\left(1-\frac{1}{2^{1+\ell-k+p}}-\frac{j}{2^{p}}\right)\Bigr]
+
|\alpha_{j}^{(p)}|\widehat{\psi}\left(1-\frac{1}{2^{\ell-k+p}}-\frac{j}{2^{p}}\right),\,\,\,\, &j&=0,1,\ldots,2^{p}-1;\\
|\alpha_{j}^{(p+1)}|&\le|\alpha_{j}^{(p)}|\Bigl[1-\widehat{\psi}\left(\frac{1}{2^{\ell-k+p}}+\frac{j}{2^{p}}\right)\Bigr]+
|\overline{\alpha}_{j}^{(p)}|\widehat{\psi}\left(\frac{1}{2^{1+\ell-k+p}}+\frac{j}{2^{p}}\right),\,\,\,\, &j&=0,1,\ldots,2^{p}-1;\\
|\alpha_{2^{p}+j}^{(p+1)}|&\le(|\overline{\alpha}_{j}^{(p)}|+|\alpha_{j}^{(p)}|)\widehat{\psi}\left(1+\frac{1}{2^{1+\ell-k+p}}+\frac{j}{2^{p}}\right),\,\,\,\, &j&=0,1,\ldots,2^{p}-1.
\end{aligned}
}

We note that these bounds depend only on $\ell-k$ and so we can iterate these inequalities from the initial (level $p=1$) bounds on $|\alpha_{0}^{(1)}|, |\overline{\alpha}_{0}^{(1)}|, |\alpha_{1}^{(1)}| $ and $|\overline{\alpha}_{1}^{(1)}|$ to compute numerical bounds on $|\alpha_{j}^{(p)}|$ and $|\overline{\alpha}_{j}^{(p)}|$, $j=0,1,\ldots,2^{p}-1.$ We can then use these bounds to compute majorants for
\beqn \label{numbnd}
\| T_{\frac{1}{2^{\ell}},p} c_m \|_\infty & \le & \sum_{j=0}^{2^p-1} |\overline{\alpha}_{j}^{(p)}|+  |\alpha_{j}^{(p)}|.
\eeqn

In order to make the development more presentable let us define $\gamma_0=1$, let $\gamma_p=(1-\hpsi(2^{-p}))$ for $p\ge1$ and define

\EQ{gamp}{
\Gamma_p = \left \{ \begin{array}{cc}
\prod_{i=1}^p \gamma_i & p \ge 1,\\
1, & p=0, \\
0, & p<0.
\end{array} \right.
}
The first $10$ values of this sequence are displayed in Table \ref{tabgauss}.

\begin{table}[h]
\centering
\begin{tabular}{|c|c|c|} \hline
Level $p$  & $\gamma_p$ & $\Gamma_p$ \\ \hline \hline
1 &  9.9(-1) & 9.9(-1) \\ \hline
2 &  7.1(-1) & 7.0(-1)  \\ \hline
3 &  2.7(-1) & 1.9(-1) \\ \hline
4 &  7.4(-2) & 1.4(-2) \\ \hline
5 &  1.9(-2)& 2.6(-4) \\ \hline
6 &  4.8(-3)& 1.2(-6) \\ \hline
7 &  1.2(-3) & 1.5(-9) \\ \hline
8 &  3.0(-4)& 4.6(-13) \\ \hline
9 &  7.5(-5) & 3.5(-17) \\ \hline
10 &  1.9(-5) & 6.9(-22) \\ \hline
\end{tabular}
\caption{Coefficients $\gamma_p$ and $\Gamma_p$.}
 \label{tabgauss}
\end{table}
Let us consider the case where $\ell -k=1,$ i.e., $\frac{1}{4}\le \frac{m}{2^{\ell}}\le \frac{1}{2},$ then appealing to (\ref{startcoefs}) we have that
\EQ{initlk1}{
\begin{aligned}
|\overline{\alpha}_{0}^{(1)}| \le \widehat{\psi}\left(\frac{3}{2}\right)\quad
 |\overline{\alpha}_{1}^{(1)}| \le \widehat{\psi}\left(\frac{1}{2}\right)\quad
|\alpha_{0}^{(1)}| \le \left(1-\widehat{\psi}\left(\frac{1}{2}\right)\right)=\Gamma_1\quad {\rm{and}}\quad |\alpha_{1}^{(1)}| \le \widehat{\psi}\left(\frac{5}{4}\right)
\end{aligned}
}
We recall from the example of approximating the constant that, in a typical error expression, the coefficients multiplying higher frequency cosines are much smaller in magnitude in relation to those multiplying the lower frequency cosines. This is evident in the above since $\widehat{\psi}\left(\frac{3}{2}\right)<<\hpsi\left(\frac{5}{4}\right) \approx 10^{-14}$  which is small in comparison to $\hpsi\left(\frac{1}{2}\right)\approx 0.007191$ and $\Gamma_{1}\approx 0.99281.$
Our aim is to develop upper bounds that are accurate to two significant figures and
so can write
\beqn \label{p1bnd}
\; \| T_{\frac{1}{2^{\ell}},1} c_m \|_{\infty}  \le |\alpha_{0}^{(1)}|+ |\overline{\alpha}_{1}^{(1)}| \le \Gamma_1+\hpsi\left(\frac{1}{2}\right).
\eeqn
If we now use (\ref{boundedcoeffs}), with $\ell-k=1$ we have the following bounds on the $p=2$ coefficients:
\EQ{initlk2-first}{
\begin{aligned}
\Bigl\{ |\overline \alpha_{0}^{(2)}|  \le & \hpsi\left(\frac{7}{4}\right) |\alpha_{0}^{(1)}| \le \hpsi\left(\frac{7}{4}\right) \Gamma_1 \Bigr\}&\quad  |\alpha_{0}^{(2)}|  \le & \left(1-\hpsi\left(\frac{1}{4}\right)\right)|\alpha_{0}^{(1)}| \le \Gamma_2 &\\
  \Bigl\{|\overline \alpha_{1}^{(2)}| \le &  \hpsi\left(\frac{5}{4}\right)|\overline \alpha_{1}^{(1)}| \le \hpsi\left(\frac{5}{4}\right)\hpsi\left(\frac{1}{2}\right) \Bigr\}&\quad |\alpha_{1}^{(2)}|  \le & \hpsi\left(\frac{5}{8}\right) |\overline \alpha_{1}^{(1)}|  \le \hpsi\left(\frac{1}{2}\right) \hpsi\left(\frac{5}{8}\right)&\\
   |\overline \alpha_{2}^{(2)}|  \le &  \hpsi\left(\frac{3}{4}\right) |\alpha_{0}^{(1)}| \le \hpsi\left(\frac{3}{4}\right) \Gamma_1 &\quad \Bigl\{ |\alpha_{2}^{(2)}|  \le & \hpsi\left(\frac{9}{8}\right) | \alpha_{0}^{(1)}|  \le \Gamma_{1} \hpsi\left(\frac{9}{8}\right)\Bigr\}&\\
 |\overline \alpha_{3}^{(2)}|  \le  &\left(1-\hpsi\left(\frac{3}{8}\right)\right) |\overline \alpha_{1}^{(1)}| \le  \Gamma_1\hpsi\left(\frac{1}{2}\right) &\quad
\Bigl\{|\alpha_{3}^{(2)}|  \le & \hpsi\left(\frac{13}{8}\right)  |\overline \alpha_{1}^{(1)}|  \le \hpsi\left(\frac{1}{2}\right) \hpsi\left(\frac{13}{8}\right)\Bigr\}.&
\end{aligned}
}
The inequalities in the braces are for those coefficients that multiply the higher frequency  cosines and, as can easily be verified, their magnitudes are negligible in comparison to the remaining coefficients. Thus, given that we seek bounds that are correct to 2 significant figures, we can write
\beqn \label{p2bnd}
\begin{aligned}
\| T_{\frac{1}{2^{\ell}},2} c_m \|_{\infty}  & \le |\alpha_{0}^{(2)}|  +  |\overline \alpha_{3}^{(2)}| + |\overline \alpha_{2}^{(2)}|+ |\alpha_{0}^{(2)}|  \le  \Gamma_2+\hpsi\left(\frac{1}{2}\right)\Gamma_1+\Gamma_1 \hpsi\left(\frac{3}{4}\right)+\hpsi\left(\frac{1}{2}\right)\hpsi\left(\frac{5}{8}\right)\\&\le \Gamma_2+\hpsi\left(\frac{1}{2}\right)\Gamma_1+2\Gamma_1 \hpsi\left(\frac{3}{4}\right).
\end{aligned}
\eeqn
where, in the final line, we have used that $\Gamma_1 \hpsi\left(\frac{3}{4}\right)>\hpsi\left(\frac{1}{2}\right)\hpsi\left(\frac{5}{8}\right).$


One more iteration is sufficient to observe the pattern of coefficients generated, this time we will only consider the eight coefficients that multiply the lower frequency cosines, namely $|\overline \alpha_{4+j}^{(3)}|$ and  $|\alpha_{j}^{(3)}|$ for $j=0,1,2,3.$
\[
\begin{aligned}
\Bigl\{ |\overline \alpha_{4}^{(3)}|  \le & \hpsi\left(\frac{7}{8}\right) |\alpha_{0}^{(3)}| \le \Gamma_{2}\hpsi\left(\frac{7}{8}\right)\Bigr\} &\\
\Bigl\{  |\overline \alpha_{5}^{(3)}| \le &  \hpsi\left(\frac{5}{8}\right)|\alpha_{1}^{(2)}| \le \hpsi^{2}\left(\frac{5}{8}\right)\hpsi\left(\frac{1}{2}\right) \Bigr\}&\\
   \Bigl[|\overline \alpha_{6}^{(3)}|  \le &  \left(1-\hpsi\left(\frac{3}{8}\right)\right) |\overline\alpha_{2}^{(2)}| \le  \Gamma_1\hpsi\left(\frac{3}{4}\right) \Gamma_1\Bigr] &\\
 |\overline \alpha_{7}^{(2)}|  \le  &\left(1-\hpsi\left(\frac{3}{16}\right)\right) |\overline \alpha_{3}^{(2)}| \le \left(1-\hpsi\left(\frac{3}{16}\right)\right) \Gamma_1\hpsi\left(\frac{1}{2}\right)\le \hpsi\left(\frac{1}{2}\right)\Gamma_{2} &\\
  |\alpha_{0}^{(3)}|  \le & \left(1-\hpsi\left(\frac{1}{8}\right)\right)|\alpha_{0}^{(2)}| \le \Gamma_3 &\\
   \Bigl[|\alpha_{1}^{(3)}|  \le & \left(1-\hpsi\left(\frac{5}{8}\right)\right) |\alpha_{1}^{(2)}|  \le \left(1-\hpsi\left(\frac{5}{8}\right)\right)\hpsi\left(\frac{1}{2}\right) \hpsi\left(\frac{5}{8}\right)\Bigr]&\\
    \Bigl\{|\alpha_{2}^{(3)}|  \le & \hpsi\left(\frac{9}{16}\right) | \overline \alpha_{2}^{(2)}|  \le \hpsi\left(\frac{9}{16}\right)\Gamma_{1}\hpsi\left(\frac{3}{4}\right)\Bigr\}&\\
\Bigl\{|\alpha_{3}^{(2)}|  \le & \hpsi\left(\frac{13}{16}\right)  |\overline \alpha_{3}^{(2)}|  \le \Gamma_{1}\hpsi\left(\frac{1}{2}\right)
\hpsi\left(\frac{13}{16}\right)\Bigr\}.&
\end{aligned}
\]

It is easy to verify that, of the two inequalities displayed in square brackets above, the largest upper bound is $\Gamma_1\hpsi\left(\frac{3}{4}\right) \Gamma_1$, similarly for the four displayed in braces, the largest upper bound is $\Gamma_{2}\hpsi\left(\frac{7}{8}\right).$ This allows us to conclude that
\[ \label{p3bnd}
\| T_{\frac{1}{2^{\ell}},3} c_m \|_{\infty} \le \Gamma_3+\Gamma_0\hpsi\left(\frac{1}{2}\right)\Gamma_2+ 2\Gamma_1 \hpsi\left(\frac{3}{4}\right) \Gamma_1 + 4\Gamma_2 \hpsi\left(\frac{7}{8}\right)\Gamma_{0} .
\]
Iterating the above procedure we obtain the following bound for $T_{\frac{1}{2^{\ell}},p} c_m$ when $\ell-k=1$:
\beqn \label{theorlk1}
\| T_{\frac{1}{2^{\ell}},p} c_m \|_{\infty}  \le \mu_{1,p}:=\Gamma_{p} + \sum_{j=0}^{p-1} 2^{j}\Gamma_{j} \hpsi\left(1-\frac{1}{2^{j+1}}\right) \Gamma_{p-1-j}.
\eeqn

We can perform the same analysis for $\ell-k =n,$ i.e., when we have  $\frac{1}{2^{n+1}}\le \frac{m}{2^{\ell}}< \frac{1}{2^{n}},$ and we obtain more generally
\beqn \label{theorup-exact}
\| T_{\frac{1}{2^{\ell}},p} c_m \|_{\infty}  & \le & \mu_{n,p}:= \frac{\Gamma_{n+p-1}}{\Gamma_{n-1}} + \sum_{j=0}^{p-1} 2^{j} \frac{\Gamma_{n-1+j}}{\Gamma_{n-1}}\hpsi\left(1-\frac{1}{2^{j+n}}\right) \Gamma_{p-1-j}.
\eeqn



In Table \ref{compare} we compute the bound above with the numerical computation of $\| T_{\frac{1}{2^{\ell}},p} c_m \|_{\infty}$   and observe that these are accurate estimates. Indeed, the results show that, to an accuracy of $\sim10^{-15},$ we can closely bound the errors in the multilevel iteration given in Proposition~\ref{propiter} when the frequency is less than half of the sampling rate. With some elementary manipulation one can easily show that
\[
\mu_{n,p}\le \gamma_{n-1}\mu_{1,p}, \quad p\ge 1, \,\, {\rm{and}}\,\,\, n\ge1,
\]
where equality holds for $n=1.$ This observation allows the error bounds to be expressed in a simpler form, i.e., solely in terms of an appropriate multiple of the sequence $ \mu_{1,p},$ and it also
 motivates the following assumption concerning a typical iteration:

\begin{assumption} \label{ass1}
Let $\ell \in \NN$ and  suppose $2^{k-1} <m \le 2^k \le 2^{\ell-1}$, so that $\frac{1}{2^{\ell-k}}\le \frac{m}{2^{\ell-1}}\le \frac{1}{2^{\ell-k-1}}.$ Then there exists a sequence of positive decreasing numbers $\mu_1>\mu_2>\cdots$ such that
\beqns
\| T_{\frac{1}{2^{\ell}},p} c_m \|_{\infty}  & \le & \gamma_{\ell-k-1} \mu_p, \quad p\ge 1.
\eeqns
\end{assumption}

\begin{table}[h]
\centering
\caption{Comparison of numerical and analytic bounds}
\label{numancom}
\begin{tabular}{|c||c|c|c|c|c|c|} \hline
$p$ & $\ell-k=1$ & $\ell-k=1$ & $\ell-k=2$ & $\ell-k=2$ & $\ell-k=3$ & $\ell-k=3$ \\
 & analytical $(\mu_{1,p})$ & numerical & analytical  $(\mu_{2,p})$ & numerical  & analytical $(\mu_{3,p})$ & numerical \\ \hline \hline
1 & 1.0 & 1.0 & 7.1(-1) & 7.1(-1)& 2.7(-1) & 2.7(-1) \\ \hline
2 & 7.1(-1) & 7.1(-1) & 1.9(-1) & 1.9(-1)& 2.0(-2) & 2.0(-2)\\ \hline
3 & 1.9(-1)& 1.9 (-1) & 1.4(-2) & 1.4(-2)& 3.8(-4) &  3.7(-4) \\ \hline
4 & 1.5(-2) & 1.4(-2) & 2.7(-4) & 2.7(-4)& 1.9(-6) & 1.9(-6)  \\ \hline
5 & 3.7(-4) & 2.7(-4) & 1.6(-6) & 1.5(-6)& 9.4(-9) & 7.3(-9)   \\ \hline
6 & 3.8(-6) & 1.3(-6) & 1.6(-8) & 7.9(-9)& 4.3(-10) & 1.9(-10) \\ \hline
7 & 3.9(-8) & 2.2(-9) & 6.3(-10) & 1.2(-10)& 1.7(-11) & 4.4(-12)  \\ \hline
8 & 1.3(-9) & 4.5(-11) & 2.4(-11) & 3.1(-12)& 4.3(-13) & 8.6(-14) \\ \hline
9 & 4.8(-11) & 1.9(-12) & 6.1(-13) & 6.0(-14)& 6.4(-15) & 1.0(-15) \\ \hline
10 & 1.2(-12) & 4.8(-14) & 9.0(-15) & 7.4(-16)& &   \\ \hline
11 & 1.8(-14) & 6.2(-16) &  & & &   \\ \hline
\end{tabular}
 \label{compare}
\end{table}





\section{General convergence of the discrete algorithm} \label{fullalg}

In this section we will examine the behaviour of the multilevel method applied to $c_{m}$ if algorithm is started by sampling at the integers, i.e., the initial sample rate $=1.$ We begin by considering the how the the error rates develop as the sample rates grow beyond $m$ the size of the cosine frequency. In order to motivate the general analysis it is instructive to examine the performance on a concrete example.\\

$\diamond$ {{\bf{The performance of the multilevel method applied to $c_{7}$}}.} \\

We begin by recalling the definition of the multilevel error operator as the algorithm proceeds,
\[
\begin{aligned}
{\rm{Initialisation:}} \,\,{\rm{Identity}}\,\,{\rm{operator:}}\quad & M_{h,0}(f) = f , \quad h>0.\\
{\rm{Level}} \,\,{\rm{one}}\,\,{\rm{error:}}\quad & M_{\frac{1}{n},1}(f) = E_{\frac{1}{n}}(f)=Q_{\frac{1}{n}} f-f.\\
{\rm{Level}} \,\,p\,\, (\ge2)\,\,{\rm{error:}}\quad & M_{\frac{1}{n},p}(f) = E_{\frac{1}{n2^{p-1}}}(M_{\frac{1}{n},p-1}(f)).
\end{aligned}
\]
Considering the case where $f=c_{7}$ we have\\

$\diamond$ {{\bf{Level one error}}.}
\[
M_{1,1}c_{7}=Q_{1}c_{7}-c_{7}=Q_{1}c_{7\md 1}-c_{7}=M_{\frac{1}{2},0} Q_{1}c_{7\md 1}-c_{7},
\]
where we have used the aliasing result (\ref{mods}) and have also introduced $M_{\frac{1}{2},0}$ to represent the identity operator, this form will be helpful in developing the further error expressions as the algorithm progresses.\\

$\diamond$ {{\bf{Level two error}}.} Using the definition of the multilevel error, together with the  aliasing result (\ref{mods}) we have
\[
\begin{aligned}
M_{1,2}c_{7}=E_{\frac{1}{2}}M_{\frac{1}{2},0} Q_{1}c_{7\md 1}-E_{\frac{1}{2}}c_{7}&=M_{\frac{1}{2},1}Q_{1}c_{7\md 1}-Q_{\frac{1}{2}}c_{7}+c_{7}\\
&=M_{\frac{1}{2},1}Q_{1}c_{7\md 1}-M_{\frac{1}{4},0}Q_{\frac{1}{2}}c_{7\md 2} +c_{7},
\end{aligned}
\]
where, in this case, $M_{\frac{1}{4},0}$ is introduced to act as the identity operator.\\

 $\diamond$ {{\bf{Level three error}}.} In the same fashion as before we have
\[
\begin{aligned}
M_{1,3}c_{7}&=E_{\frac{1}{4}}M_{\frac{1}{2},1}Q_{1}c_{7\md 1}-E_{\frac{1}{4}}M_{\frac{1}{4},0}Q_{\frac{1}{2}}c_{7\md 2} +E_{\frac{1}{4}}c_{7}\\
&=M_{\frac{1}{2},2}Q_{1}c_{7\md 1}-M_{\frac{1}{4},1}Q_{\frac{1}{2}}c_{7\md 2} +Q_{\frac{1}{4}}c_{7}-c_{7}\\
&=M_{\frac{1}{2},2}Q_{1}c_{7\md 1}-M_{\frac{1}{4},1}Q_{\frac{1}{2}}c_{7\md 2} +M_{\frac{1}{8},0}Q_{\frac{1}{4}}c_{7\md 4}-c_{7},
\end{aligned}
\]
where, in this case, $M_{\frac{1}{8},0}$ is introduced to act as the identity operator.\\

 $\diamond$ {{\bf{Level four error}}.} At this iteration we have
\[
\begin{aligned}
M_{1,4}c_{7}&=E_{\frac{1}{8}}M_{\frac{1}{2},2}Q_{1}c_{7\md 1}-E_{\frac{1}{8}}M_{\frac{1}{4},1}Q_{\frac{1}{2}}c_{7\md 2} +E_{\frac{1}{8}}M_{\frac{1}{8},0}Q_{\frac{1}{4}}c_{7\md 4}-E_{\frac{1}{8}}c_{7},\\
&=M_{\frac{1}{2},3}Q_{1}c_{7\md 1}-M_{\frac{1}{4},2}Q_{\frac{1}{2}}c_{7\md 2} +M_{\frac{1}{8},1}Q_{\frac{1}{2}}c_{7\md 4} -
Q_{\frac{1}{8}}c_{7}+c_{7}\\
&=M_{\frac{1}{2},3}Q_{1}c_{7\md 1}-M_{\frac{1}{4},2}Q_{\frac{1}{2}}c_{7\md 2} +M_{\frac{1}{8},1}Q_{\frac{1}{2}}c_{7\md 4} -
M_{\frac{1}{16},0}Q_{\frac{1}{8}}c_{7\md8}+c_{7}\,
\end{aligned}
\]

 $\diamond$ {{\bf{Level five error}}.} At this iteration we have
\[
\begin{aligned}
M_{1,5}c_{7}&=E_{\frac{1}{16}}M_{\frac{1}{2},3}Q_{1}c_{7\md 1}-E_{\frac{1}{16}}M_{\frac{1}{4},2}Q_{\frac{1}{2}}c_{7\md 2} +E_{\frac{1}{16}}M_{\frac{1}{8},1}Q_{\frac{1}{2}}c_{7\md 4} -
E_{\frac{1}{16}}M_{\frac{1}{16},0}Q_{\frac{1}{8}}c_{7\md8}+E_{\frac{1}{16}}c_{7},\\
&=M_{\frac{1}{2},4}Q_{1}c_{7\md 1}-M_{\frac{1}{4},3}Q_{\frac{1}{2}}c_{7\md 2} +M_{\frac{1}{8},2}Q_{\frac{1}{2}}c_{7\md 4} -M_{\frac{1}{16},1}Q_{\frac{1}{8}}c_{7\md8}+
\underbrace{Q_{\frac{1}{16}}c_{7}-c_{7}}_{=M_{\frac{1}{16},1}c_{7}}\\
\end{aligned}
\]

We notice that at the fifth level of iteration the cosine frequency ($=7$) is now less that half of the current sample rate $(=2^{4}=16$), Thus, we can appeal to the analysis of the previous section to track how the algorithm evolves for this part of the error expression. Using this example as a template it is straightforward to establish the following more general result.

\begin{proposition} \label{highfreq}
Let $\ell$  and $m$ be non-negative integers that satisfy $2^{\ell}\le m <2^{\ell+1}.$ Then,  we have
\beqns
M_{1,p} c_m & = & \sum_{j=0}^{p-1}(-1)^{j} M_{\frac{1}{2^{j+1}},p-(j+1)} Q_{\frac{1}{2^{j}}} c_{m (\md 2^j)}+(-1)^{p}c_{m}\quad {\rm{for}}\,\,1 \le p \le \ell+2
\eeqns
and
\beqns
M_{1,p} c_m & = & \sum_{j=0}^{\ell+1}(-1)^{j} M_{\frac{1}{2^{j+1}},p-(j+1) }Q_{\frac{1}{2^{j}}} c_{m (\md 2^j)} + (-1)^{\ell+2}M_{\frac{1}{2^{\ell+2}},p-(\ell+2)}c_{m} \quad {\rm{for}}\,\,p\ge \ell+3.
\eeqns
\end{proposition}

In order to deliver error bounds for the full multilevel algorithm we need to track the build up of terms, captured in Proposition \ref{highfreq}, that arise from the early iterations where the sample rate is small in comparison to the cosine frequency. To shed more light onto this we will continue to develop the application to $c_{7}$  as a prototype. For this example we will take the initial sample rate to be $8=2^{3},$ and we will derive an expression  for the term $ M_{\frac{1}{16},p}Q_{\frac{1}{8}}c_{7}.$ In what follows we will make heavy use of  identity (\ref{Qothers}). We have already observed, with the example of the constant, how the multilevel error expansions grow rapidly as the algorithm progresses and, in order to exert some control on the development we will aim monitor our example  $M_{\frac{1}{16},p}Q_{\frac{1}{8}}c_{7}$ by writing it in terms of  multilevel expressions applied to those cosines whose frequencies are less than half of the current sample rate (for such terms we can appeal to Proposition \ref{propiter}) while monitoring the remainder terms. As we shall see, the coefficients that multiply the incoming cosine terms (from one level to the next) decay very quickly and, in fact, only  three iterations are required before  the magnitude  of these coefficients drops below $\epsilon.$\\

$\diamond$ {{\bf{Level one error}}.}
\[
\begin{aligned}
&M_{\frac{1}{16},1}Q_{\frac{1}{8}}c_{7}=Q_{\frac{1}{16}}Q_{\frac{1}{8}}c_{7}-Q_{\frac{1}{8}}c_{7}\\
&=\widehat{\psi}\left(\frac{1}{8}\right)\left(Q_{\frac{1}{16}}c_{1}-c_{1}\right)
+\widehat{\psi}\left(\frac{7}{8}\right)\left(Q_{\frac{1}{16}}c_{7}-c_{7}\right)\\
&+\widehat{\psi}\left(\frac{9}{8}\right)\left(Q_{\frac{1}{16}}c_{9}-c_{9}\right)
+\widehat{\psi}\left(\frac{15}{8}\right)\left(Q_{\frac{1}{16}}c_{15}-c_{15}\right)+Q_{\frac{1}{16}}g_{0}-g_{0}\\
&=\widehat{\psi}\left(\frac{1}{8}\right)M_{\frac{1}{16},1}c_{1}+\widehat{\psi}\left(\frac{7}{8}\right)M_{\frac{1}{16},1}c_{7}\\
&+\widehat{\psi}\left(\frac{9}{8}\right)\Bigl[\widehat{\psi}\left(\frac{7}{16}\right)c_{7}
+\left(\widehat{\psi}\left(\frac{9}{16}\right)-1\right)c_{9}\Bigr]+\widehat{\psi}\left(\frac{9}{8}\right)\Bigl\{\widehat{\psi}\left(\frac{23}{16}\right)c_{23}+
\widehat{\psi}\left(\frac{25}{16}\right)c_{25}+g_{1}^{(1)}\Bigr\}\\
&+\widehat{\psi}\left(\frac{15}{8}\right)\Bigl[\widehat{\psi}\left(\frac{1}{16}\right)c_{1}
+\left(\widehat{\psi}\left(\frac{15}{16}\right)-1\right)c_{15}\Bigr]+\widehat{\psi}\left(\frac{15}{8}\right)\Bigl\{\widehat{\psi}\left(\frac{17}{16}\right)c_{17}+
\widehat{\psi}\left(\frac{31}{16}\right)c_{31}+g_{2}^{(1)}\Bigr\}\\
&+Q_{\frac{1}{16}}g_{0}-g_{0}\\
\end{aligned}
\]
where \[
g_{1}=\widehat{\psi}\left(\frac{9}{8}\right)g_{1}^{(1)}+\widehat{\psi}\left(\frac{15}{8}\right)g_{1}^{(1)}+Q_{\frac{1}{16}}g_{0}-g_{0}
\]
for which we have
\[
\|g_{1}\|_{\infty}\le 2\widehat{\psi}\left(1\right)\epsilon+A\epsilon \le 2A\epsilon.
\]

Observing this first iteration we notice that of the $8$ newly introduced cosine terms the $4$ lower frequency examples (those below 16 and appearing in square brackets) are  multiplied by coefficients that are bounded by $\widehat{\psi}\left(1\right),$ whereas  the remaining $4$ higher frequency cosines (those above 16 and appearing in braces)  are  multiplied by coefficients that are bounded by   $\widehat{\psi}\left(1\right)^{2}.$  This allows us to write
\EQ{nice1}{
\begin{aligned}
M_{\frac{1}{16},1}Q_{\frac{1}{8}}c_{7}=\rho_{-1}^{(0)}M_{\frac{1}{16},1}c_{1}+\rho_{0}^{(0)}M_{\frac{1}{16},1}c_{7}
&+\rho_{-2}^{(1)}c_{9}+\rho_{-1}^{(1)}c_{1}+\rho_{0}^{(1)}c_{7}+\rho_{1}^{(1)}c_{15}\\
&+\gamma_{-4}^{(2)}c_{25}+\gamma_{-3}^{(2)}c_{17}+\gamma_{2}^{(2)}c_{23}+\gamma_{3}^{(2)}c_{31}+g_{1},
\end{aligned}
}
where $|\rho_{i}^{(0)}|<1$ ($i=-1,0$), $|\rho_{i}^{(1)}|<\widehat{\psi}\left(1\right)$ ($i=-2,-1,0,1$), and $|\gamma_{i}^{(2)}|<\widehat{\psi}\left(1\right)^{2}$ ($i=-4,-3,2,3$). Using (\ref{nice1}) as base for progressing we move onto the next level.\\

$\diamond$ {{\bf{Level two error}}.}

\[
\begin{aligned}
M_{\frac{1}{16},2}Q_{\frac{1}{8}}c_{7}&=\rho_{-1}^{(0)}M_{\frac{1}{16},2}c_{1}+\rho_{0}^{(0)}M_{\frac{1}{16},2}c_{7}\\
&+\rho_{-2}^{(1)}M_{\frac{1}{32},1}c_{9}+\rho_{-1}^{(1)}M_{\frac{1}{32},1}c_{1}+\rho_{0}^{(1)}M_{\frac{1}{32},1}c_{7}+\rho_{1}^{(1)}M_{\frac{1}{32},1}c_{15}\\
&+\gamma_{-4}^{(2)}\Bigl[\widehat{\psi}\left(\frac{7}{32}\right)c_{7}+\left(\widehat{\psi}\left(\frac{25}{32}\right)-1\right)c_{25}\Bigr]+
\gamma_{-4}^{(2)}\Bigl\{\widehat{\psi}\left(\frac{39}{32}\right)c_{39}+\widehat{\psi}\left(\frac{57}{32}\right)c_{49}\Bigr\}\\
&+\gamma_{-3}^{(2)}\Bigl[\widehat{\psi}\left(\frac{15}{32}\right)c_{15}+\left(\widehat{\psi}\left(\frac{17}{32}\right)-1\right)c_{17}\Bigr]+
\gamma_{-3}^{(2)}\Bigl\{\widehat{\psi}\left(\frac{47}{32}\right)c_{47}+\widehat{\psi}\left(\frac{49}{32}\right)c_{49}\Bigr\}\\
&+\gamma_{2}^{(2)}\Bigl[\widehat{\psi}\left(\frac{9}{32}\right)c_{9}+\left(\widehat{\psi}\left(\frac{23}{32}\right)-1\right)c_{23}\Bigr]+
\gamma_{2}^{(2)}\Bigl\{\widehat{\psi}\left(\frac{41}{32}\right)c_{41}+\widehat{\psi}\left(\frac{55}{32}\right)c_{55}\Bigr\}\\
&+\gamma_{3}^{(2)}\Bigl[\widehat{\psi}\left(\frac{1}{32}\right)c_{1}+\left(\widehat{\psi}\left(\frac{31}{32}\right)-1\right)c_{31}\Bigr]+
\gamma_{3}^{(2)}\Bigl\{\widehat{\psi}\left(\frac{33}{32}\right)c_{33}+\widehat{\psi}\left(\frac{63}{32}\right)c_{63}\Bigr\}+g_{2}
\end{aligned}
\]
where
\[
g_{2}=\gamma_{-4}^{(2)}g_{-4}^{(2)}+\gamma_{-3}^{(2)}g_{-3}^{(2)}+\gamma_{2}^{(2)}g_{2}^{(2)}+\gamma_{3}^{(2)}g_{3}^{(2)}+Q_{\frac{1}{32}}g_{1}-g_{1}.
\]
for which we have
\[
\|g_{2}\|_{\infty}\le 4\widehat{\psi}\left(1\right)^{2}\epsilon+A\|g_{1}\|_{\infty}\le (4\widehat{\psi}\left(1\right)^{2}+2A\widehat{\psi}\left(1\right)\epsilon+A^{2})\epsilon=
 \le 2A^{2}\epsilon.
\]
Just as in the first level we notice that of the  16 newly introduced  cosine terms, the half in the lower frequency set (those below 32 and in square brackets) are  multiplied by coefficients  of magnitude less than $\widehat{\psi}\left(1\right)^{2},$ whereas the remaining higher frequency cosines (those above 32 and in the braces)  are  multiplied by coefficients  of magnitude less than $\widehat{\psi}\left(1\right)^{3}.$  In view of this we can express this iteration as
\EQ{nice2}{
\begin{aligned}
M_{\frac{1}{16},2}Q_{\frac{1}{8}}c_{7}&=\rho_{-1}^{(0)}M_{\frac{1}{16},2}c_{1}+\rho_{0}^{(0)}M_{\frac{1}{16},2}c_{7}\\
&+\rho_{-2}^{(1)}M_{\frac{1}{32},1}c_{9}+\rho_{-1}^{(1)}M_{\frac{1}{32},1}c_{1}+\rho_{0}^{(1)}M_{\frac{1}{32},1}c_{7}+\rho_{1}^{(1)}M_{\frac{1}{32},1}c_{15}\\&+\rho_{-4}^{(2)}c_{25}+\rho_{-3}^{(2)}c_{17}+\rho_{-2}^{(2)}c_{9}+\rho_{-1}^{(2)}c_{1}+\rho_{0}^{(2)}c_{7}+\rho_{1}^{(2)}c_{15}+\rho_{2}^{(2)}c_{23}+\rho_{3}^{(2)}c_{31}\\
&+\gamma_{-8}^{(3)}c_{57}+\gamma_{-7}^{(3)}c_{49}+\gamma_{-6}^{(3)}c_{41}+\gamma_{-5}^{(3)}c_{33}+\gamma_{4}^{(3)}c_{39}+\gamma_{5}^{(3)}c_{47}+\gamma_{6}^{(3)}c_{55}+\gamma_{7}^{(3)}c_{63}+g_{2},
\end{aligned}
}
where, as before, where $|\rho_{i}^{(0)}|<1$ ($i=-1,0$), $|\rho_{i}^{(1)}|<\widehat{\psi}\left(1\right)$ ($i=-2,-1,0,1$), but now we have, in addition, that
$|\rho_{i}^{(2)}|<\widehat{\psi}\left(1\right)^{2}$ ($i=-4,\ldots,3$), and  $|\gamma_{-8+i}^{(3)}|<\widehat{\psi}\left(1\right)^{3},$ $|\gamma_{7-i}^{(3)}|<\widehat{\psi}\left(1\right)^{3},$  for ($i=0,1,2,3$).
With this formula established we can move to the next level.\\

$\diamond$ {{\bf{Level three error}}.}

\[
\begin{aligned}
M_{\frac{1}{16},3}Q_{\frac{1}{8}}c_{7}&=\rho_{-1}^{(0)}M_{\frac{1}{16},3}c_{1}+\rho_{0}^{(0)}M_{\frac{1}{16},3}c_{7}\\
&+\rho_{-2}^{(1)}M_{\frac{1}{32},2}c_{9}+\rho_{-1}^{(1)}M_{\frac{1}{32},2}c_{1}+\rho_{0}^{(1)}M_{\frac{1}{32},2}c_{7}+\rho_{1}^{(1)}M_{\frac{1}{32},2}c_{15}\\
&+\rho_{-4}^{(2)}M_{\frac{1}{64},1}c_{25}+\rho_{-3}^{(2)}M_{\frac{1}{64},1}c_{17}+\rho_{-2}^{(2)}M_{\frac{1}{64},1}c_{9}+\rho_{-1}^{(2)}M_{\frac{1}{64},1}c_{1}\\&+\rho_{0}^{(2)}M_{\frac{1}{64},1}c_{7}+\rho_{1}^{(2)}M_{\frac{1}{64},1}c_{15}+\rho_{2}^{(2)}M_{\frac{1}{64},1}c_{23}+\rho_{3}^{(2)}M_{\frac{1}{64},1}c_{31}\\
&+\gamma_{-8}^{(3)}\Bigl[\widehat{\psi}\left(\frac{7}{64}\right)c_{7}+\left(\widehat{\psi}\left(\frac{57}{64}\right)-1\right)c_{57}\Bigr]+
\gamma_{-8}^{(3)}\Bigl\{\widehat{\psi}\left(\frac{71}{64}\right)c_{71}+\widehat{\psi}\left(\frac{121}{64}\right)c_{121}\Bigr\}\\
&+\gamma_{-7}^{(3)}\Bigl[\widehat{\psi}\left(\frac{15}{64}\right)c_{15}+\left(\widehat{\psi}\left(\frac{49}{64}\right)-1\right)c_{49}\Bigr]+
\gamma_{-7}^{(3)}\Bigl\{\widehat{\psi}\left(\frac{79}{64}\right)c_{79}+\widehat{\psi}\left(\frac{113}{64}\right)c_{113}\Bigr\}\\
&+\gamma_{-6}^{(3)}\Bigl[\widehat{\psi}\left(\frac{23}{64}\right)c_{23}+\left(\widehat{\psi}\left(\frac{41}{64}\right)-1\right)c_{41}\Bigr]+
\gamma_{-6}^{(3)}\Bigl\{\widehat{\psi}\left(\frac{87}{64}\right)c_{87}+\widehat{\psi}\left(\frac{105}{64}\right)c_{105}\Bigr\}\\
&+\gamma_{-5}^{(3)}\Bigl[\widehat{\psi}\left(\frac{31}{64}\right)c_{31}+\left(\widehat{\psi}\left(\frac{33}{64}\right)-1\right)c_{33}\Bigr]+
\gamma_{-5}^{(3)}\Bigl\{\widehat{\psi}\left(\frac{95}{64}\right)c_{95}+\widehat{\psi}\left(\frac{97}{64}\right)c_{97}\Bigr\}\\
&+\gamma_{4}^{(3)}\Bigl[\widehat{\psi}\left(\frac{25}{64}\right)c_{25}+\left(\widehat{\psi}\left(\frac{39}{64}\right)-1\right)c_{39}\Bigr]+
\gamma_{4}^{(3)}\Bigl\{\widehat{\psi}\left(\frac{89}{64}\right)c_{89}+\widehat{\psi}\left(\frac{103}{64}\right)c_{103}\Bigr\}\\
&+\gamma_{5}^{(3)}\Bigl[\widehat{\psi}\left(\frac{17}{64}\right)c_{17}+\left(\widehat{\psi}\left(\frac{47}{64}\right)-1\right)c_{47}\Bigr]+
\gamma_{5}^{(3)}\Bigl\{\widehat{\psi}\left(\frac{81}{64}\right)c_{81}+\widehat{\psi}\left(\frac{111}{64}\right)c_{111}\Bigr\}\\
&+\gamma_{6}^{(3)}\Bigl[\widehat{\psi}\left(\frac{9}{64}\right)c_{9}+\left(\widehat{\psi}\left(\frac{55}{64}\right)-1\right)c_{55}\Bigr]+
\gamma_{6}^{(3)}\Bigl\{\widehat{\psi}\left(\frac{73}{64}\right)c_{73}+\widehat{\psi}\left(\frac{119}{64}\right)c_{119}\Bigr\}\\
&+\gamma_{7}^{(3)}\Bigl[\widehat{\psi}\left(\frac{1}{64}\right)c_{1}+\left(\widehat{\psi}\left(\frac{63}{64}\right)-1\right)c_{63}\Bigr]+
\gamma_{7}^{(3)}\Bigl\{\widehat{\psi}\left(\frac{65}{64}\right)c_{65}+\widehat{\psi}\left(\frac{127}{64}\right)c_{127}\Bigr\}+\widetilde{g}_{3},
\end{aligned}
\]
where
\[
\widetilde{g}_{3}=\sum_{i=0}^{3}\left(\gamma_{-8+i}^{(3)}g_{-8+i}^{(3)}+\gamma_{7-i}^{(3)}g_{7-i}^{(3)}\right)+Q_{\frac{1}{64}}g_{2}-g_{2}.
\]
In our development so far we have identified the remainder terms as those whose magnitude can be bounded by a multiple of $\epsilon.$ In the above expression we continue to observe the pattern that out of the  32 newly introduced cosine terms, the  half in the lower frequency set (those below $64$ and appearing in the square brackets) are multiplied by coefficients of magnitude less than $\widehat{\psi}\left(1\right)^{3}$ whereas the remaining higher frequency cosines (those above $64$ and appearing in the braces) are multiplied by coefficients of magnitude less than  $\widehat{\psi}\left(1\right)^{4}=\widehat{\psi}\left(2\right)= \frac{1}{3}\epsilon.$ In view of this we absorb the linear combination of these higher frequency cosines into the remainder term and write
\EQ{nice3}{
\begin{aligned}
M_{\frac{1}{16},3}Q_{\frac{1}{8}}c_{7}&=\rho_{-1}^{(0)}M_{\frac{1}{16},3}c_{1}+\rho_{0}^{(0)}M_{\frac{1}{16},3}c_{7}\\
&+\rho_{-2}^{(1)}M_{\frac{1}{32},2}c_{9}+\rho_{-1}^{(1)}M_{\frac{1}{32},2}c_{1}+\rho_{0}^{(1)}M_{\frac{1}{32},2}c_{7}+\rho_{1}^{(1)}M_{\frac{1}{32},2}c_{15}\\
&+\rho_{-4}^{(2)}M_{\frac{1}{64},1}c_{25}+\rho_{-3}^{(2)}M_{\frac{1}{64},1}c_{17}+\rho_{-2}^{(2)}M_{\frac{1}{64},1}c_{9}+\rho_{-1}^{(2)}M_{\frac{1}{64},1}c_{1}\\&+\rho_{0}^{(2)}M_{\frac{1}{64},1}c_{7}+\rho_{1}^{(2)}M_{\frac{1}{64},1}c_{15}+\rho_{2}^{(2)}M_{\frac{1}{64},1}c_{23}+\rho_{3}^{(2)}M_{\frac{1}{64},1}c_{31}\\
&+\rho_{-8}^{(3)}M_{\frac{1}{128},0}c_{57}+\rho_{-7}^{(3)}M_{\frac{1}{128},0}c_{49}+\rho_{-6}^{(3)}M_{\frac{1}{128},0}c_{41}+\rho_{-5}^{(3)}M_{\frac{1}{128},0}c_{33}\\
&+\rho_{-4}^{(3)}M_{\frac{1}{128},0}c_{25}+\rho_{-3}^{(3)}M_{\frac{1}{128},0}c_{17}+\rho_{-2}^{(3)}M_{\frac{1}{128},0}c_{9}+\rho_{-1}^{(3)}M_{\frac{1}{128},0}c_{1}\\
&+\rho_{0}^{(3)}M_{\frac{1}{128},0}c_{7}+\rho_{1}^{(3)}M_{\frac{1}{128},0}c_{15}+\rho_{2}^{(3)}M_{\frac{1}{128},0}c_{23}+\rho_{3}^{(3)}M_{\frac{1}{128},0}c_{31}\\&+\rho_{4}^{(3)}M_{\frac{1}{128},0}c_{39}+\rho_{5}^{(3)}M_{\frac{1}{128},0}c_{47}+\rho_{6}^{(3)}M_{\frac{1}{128},0}c_{55}+\rho_{7}^{(3)}M_{\frac{1}{128},0}c_{63}+g_{3}
\end{aligned}
}
where, in the above, we use the fact that $M_{\frac{1}{128},0}$ is the identity. As before the coefficients appearing in this error expression satisfy
$|\rho_{i}^{(0)}|<1$ ($i=-1,0$), $|\rho_{i}^{(1)}|<\widehat{\psi}\left(1\right)$ ($i=-2,-1,0,1$), $|\rho_{i}^{(2)}|<\widehat{\psi}\left(1\right)^{2}$ ($i=-4,\ldots,3$), and now, in addition we have, $|\rho_{i}^{(3)}|<\widehat{\psi}\left(1\right)^{3}$  $(i=-8,\ldots,7).$ The remainder term is given by
\[
\begin{aligned}
g_{3}=
&\gamma_{-16}^{(4)}c_{121}+\gamma_{-15}^{(4)}c_{113}+\gamma_{-14}^{(4)}c_{105}+\gamma_{-13}^{(4)}c_{97}+\gamma_{-12}^{(4)}c_{87}+\gamma_{-11}^{(4)}c_{81}+\gamma_{-10}^{(4)}c_{73}+\gamma_{-9}^{(4)}c_{65}\\
&+\gamma_{8}^{(4)}c_{71}+\gamma_{9}^{(4)}c_{79}+\gamma_{10}^{(4)}c_{87}+\gamma_{11}^{(4)}c_{95}+\gamma_{12}^{(4)}c_{103}+\gamma_{13}^{(4)}c_{111}+\gamma_{14}^{(4)}c_{119}+\gamma_{15}^{(4)}c_{127}\\
&+\sum_{i=0}^{3}\left(\gamma_{-8+i}^{(3)}g_{-8+i}^{(3)}+\gamma_{7-i}^{(3)}g_{7-i}^{(3)}\right)+Q_{\frac{1}{64}}g_{2}-g_{2}.
\end{aligned}
\]
where  $|\gamma_{i}^{(4)}|<\widehat{\psi}\left(1\right)^{4}=\widehat{\psi}\left(2\right)=\frac{1}{3}\epsilon$  for ($i=1,\ldots,16$). The remainder term at this level
can be bounded as follows
\[
\|g_{3}\|_{\infty}\le \frac{16}{3}\epsilon+8\widehat{\psi}\left(1\right)^{3}\epsilon+A\|g_{2}\|_{\infty}\le
\left(\frac{16}{3}+8\widehat{\psi}\left(1\right)^{3}+4A\widehat{\psi}\left(1\right)^{2}+2A^{2}\widehat{\psi}\left(1\right)+A^{3}\right)\epsilon\le 2A^{3}\epsilon.
\]

From this point on further iterations will introduce no new terms into the main part of the error expression and the magnitude of the remainder term at each level will continue to be bounded by the appropriate multiple of $\epsilon.$ Using this example as a template one can use the same analysis to prove the following more general result.


\begin{proposition} \label{mqbnd}
Suppose $m < 2^\ell$, for some $\ell \in \NN$. Then  for $p \ge 3$,
\EQ{near}{
\begin{aligned}
M_{\frac{1}{2^{\ell+1}},p}  Q_{\frac{1}{2^{\ell}}} c_m& = \sum_{i=-1}^{0} \rho_{i}^{(0)}M_{\frac{1}{2^{\ell+1}},p} c_{m+i2^{\ell}}+ \sum_{i=-2}^{1} \rho_{i}^{(1)}M_{\frac{1}{2^{\ell+2}},p-1} c_{m+i2^{\ell}}\\&+ \sum_{i=-4}^{3} \rho_{i}^{(2)}M_{\frac{1}{2^{\ell+3}},p-2} c_{m+i2^{\ell}}+ \sum_{i=-8}^{7} \rho_{i}^{(3)}M_{\frac{1}{2^{\ell+4}},p-3} c_{m+i2^{\ell}}+g_{p},
\end{aligned}
}
where $|\rho_{i}^{(0)}|<1$ ($i=-1,0$), $|\rho_{i}^{(1)}|<\widehat{\psi}\left(1\right)$ ($i=-2,-1,0,1$), $|\rho_{i}^{(2)}|<\widehat{\psi}\left(1\right)^{2}$ ($i=-4,\ldots,3$), $|\rho_{i}^{(3)}|<\widehat{\psi}\left(1\right)^{3}$  $(i=-8,\ldots,7),$ and where $\|g_{p}\|_{\infty} \le 2A^{p}\epsilon.$

\end{proposition}

\begin{proof} The result can be established by following precisely the steps as we have done to derive the error expression for   the $c_{7}$ example.

\end{proof}

The multilevel terms that appear in the main result of Proposition \ref{mqbnd} share the common characteristic that they all apply  to cosines whose frequencies are less than half of the associated sample rate. An expression for such cases was formulated  in Proposition  \ref{propiter} and using this, together with
Assumption \ref{ass1}, we can develop  (\ref{near}) as follows.

\begin{corollary} \label{lowcoeffrate1}
Suppose $m \le 2^{\ell}$, for some $\ell \in \NN$. Then  for $p \ge 3$,
\beqns
\begin{aligned}
\| M_{\frac{1}{2^{\ell+1}},p}  Q_{\frac{1}{2^{\ell}}} c_m \|_{\infty}  \le  \tilde \mu_p   :=  2 & (\mu_p+\hpsi(1)\mu_{p-1}+2\hpsi(1)^2\mu_{p-2}+4\hpsi(1)^3 \mu_{p-3}) \\
+ 2&\hpsi(1) \gamma_1 (\mu_{p-1} + \hpsi(1)\mu_{p-2}  + 2 \hpsi(1)^2\mu_{p-3} ) \\
+ 2& \hpsi(1)^2 \gamma_2  (\mu_{p-2}+ \hpsi(1)\mu_{p-3}) \\+ 2 &\hpsi(1)^3 \gamma_3 \mu_{p-3} + 32pA^p  \epsilon.
\end{aligned}
\eeqns
\end{corollary}
\begin{proof} Using  (\ref{near})
 and taking into account the bounds on the sizes of the coefficients we have that
 \[
 \begin{aligned}
\|M_{\frac{1}{2^{\ell+1}},p}  Q_{\frac{1}{2^{\ell}}} c_m\|_{\infty}& \le \sum_{i=-1}^{0} \|M_{\frac{1}{2^{\ell+1}},p} c_{m+i2^{\ell}}\|_{\infty}+ \hpsi(1)\sum_{i=-2}^{1} \|M_{\frac{1}{2^{\ell+2}},p-1} c_{m+i2^{\ell}}\|_{\infty}\\&+\hpsi(1)^{2} \sum_{i=-4}^{3} \|M_{\frac{1}{2^{\ell+3}},p-2} c_{m+i2^{\ell}}\|_{\infty}+ \hpsi(1)^{3}\sum_{i=-8}^{7} \|M_{\frac{1}{2^{\ell+4}},p-3} c_{m+i2^{\ell}}\|_{\infty}+2A^p \epsilon.
\end{aligned}
\]
Noticing that the cosine frequencies appearing above are smaller in magnitude that $2^{\ell-1}$ we have, from Proposition \ref{propiter}, that
 \[
 \begin{aligned}
\|M_{\frac{1}{2^{\ell+1}},p}  Q_{\frac{1}{2^{\ell}}} c_m\|_{\infty}& \le \sum_{i=-1}^{0} \|T_{\frac{1}{2^{\ell+1}},p} c_{m+i2^{\ell}}\|_{\infty}+ \hpsi(1)\sum_{i=-2}^{1} \|T_{\frac{1}{2^{\ell+2}},p-1} c_{m+i2^{\ell}}\|_{\infty}\\&+\hpsi(1)^{2} \sum_{i=-4}^{3} \|T_{\frac{1}{2^{\ell+3}},p-2} c_{m+i2^{\ell}}\|_{\infty}+ \hpsi(1)^{3}\sum_{i=-8}^{7} \|T_{\frac{1}{2^{\ell+4}},p-3} c_{m+i2^{\ell}}\|_{\infty}+(30pA^{p-1}+2A^p) \epsilon.
\end{aligned}
\]
Inspecting the ratio of the magnitude of cosine frequency to the appropriate sample rate we observe that in the first sum this is $<\frac{1}{2}$ for both terms, in the second sum there are two terms for which it is  $<\frac{1}{2}$ and for the remaining two terms  it is $<\frac{1}{4},$ in the third sum the partition is that four terms have a ratio that is $<\frac{1}{2},$ two have a ratio that is $<\frac{1}{4}$ and two have a ratio  $\frac{1}{8}.$  The pattern continues for the fourth term where eight terms have a ratio that is $<\frac{1}{2},$  four  $<\frac{1}{4},$ two $<\frac{1}{8}$ and two $<\frac{1}{16}.$   Assumption~\ref{ass1}  states that when the ratio of cosine frequency to sample rate is  $< \frac{1}{2^{n}}$, then the level $p$ truncation term is bounded by $\gamma_{n-1} \mu_p,$ employing this gives
 \[
 \begin{aligned}
\|M_{\frac{1}{2^{\ell+1}},p}  Q_{\frac{1}{2^{\ell}}} c_m\|_{\infty}& \le 2\gamma_{0}\mu_{p}+ \hpsi(1)(2\gamma_{0}+2\gamma_{1})\mu_{p-1}+\hpsi(1)^{2} (4\gamma_{0}+2\gamma_{1}+2\gamma_{2})\mu_{p-2}\\&+ \hpsi(1)^{3}(8\gamma_{0}+4\gamma_{1}+2\gamma_{2}+2\gamma_{3})\mu_{p-3}+32pA^p \epsilon.
\end{aligned}
\]
The $\epsilon$ component of the above bound follows from the fact that $A^{p-1}<A^{p}$ and $2A^{p}<2pA^{2}$. It is straightforward to verify that, after some rearrangement of terms, and the fact that $\gamma_{0}=1,$ the above bound is precisely that quoted in the corollary.
\end{proof}

\begin{table}[h]
\centering
\label{mu}
\begin{tabular}{|c|c|c|c|c|c|c|c|c|c|c|c|} \hline
$p$ & 3 & 4 & 5 & 6 & 7 & 8& 9 & 10 & 11 & 12 \\ \hline
$\tilde \mu_p$   & 3.8(-1) & 3.0(-2)&  7.4(-4) & 7.5(-6)& 7.7(-8) & 2.5(-9) & 9.7(-11)& 2.4(-12)& 3.6(-14)& 2.9(-16)  \\ \hline
\end{tabular}
\caption{Values of $\tilde \mu_p$ based on the values of $\mu_{p,1}$ in Table~\ref{numancom}}
\end{table}

If we take the sequence $\mu_{p}$ to match  $\mu_{p,1},$ as generated in the previous section and given by (\ref{theorlk1}), then the corresponding   values of $ \tilde \mu_p $ are presented in Table  \ref{mu} and we observe that these decay quickly. We can now combine the result above with Proposition~\ref{highfreq} to yield the following general error bound.
\begin{corollary} \label{highfreqbnd}
Let $\ell \in \NN$ and $m=2^\ell+n$ for some $0 \le n < 2^{\ell}$. Then, for $p \ge \ell+5$,
$$
 \| M_{1,p} c_m \|_{\infty} \le \sigma_{p,\ell}:=\sum_{j=0}^{\ell+1} \tilde \mu_{p-j-1} +\mu_{p-\ell-2}+(p-\ell-2)A^{p-\ell-3}\epsilon.
$$
\end{corollary}
\begin{proof} From Proposition~\ref{highfreq} we have

\EQ{fullser}{
\|M_{1,p} c_m\|_{\infty}  \le  \sum_{j=0}^{\ell+1}\|M_{\frac{1}{2^{j+1}},p-(j+1) }Q_{\frac{1}{2^{j}}} c_{m (\md 2^j)}\|_{\infty} + \|M_{\frac{1}{2^{\ell+2}},p-(\ell+2)}c_{m}\|_{\infty}\quad {\rm{for}}\,\,p\ge \ell+3.
}

Since $2^{\ell} \le m < 2^{\ell+1}$ and consequently  $\frac{m}{2^{\ell+2}}<\frac{1}{2}$, we can use Proposition~\ref{propiter} together with Assumptions~\ref{ass1}, to yield
\beqn
\|M_{\frac{1}{2^{\ell+2}},p-(\ell+2)}c_{m}\|_{\infty} & \le & \gamma_0 \mu_{p-\ell-2}+(p-\ell-2)A^{p-\ell-3}\epsilon=\mu_{p-\ell-2}+(p-\ell-2)A^{p-\ell-3}\epsilon. \label{topterm}
\eeqn
Furthermore, since $p-(\ell+2)\ge 3$ we can appeal to Corollary~\ref{lowcoeffrate1}, to yield
\beqn
\| M_{\frac{1}{2^{j+1}},p-j-1} Q_{\frac{1}{2^{j}}} c_{m (\md 2^j)} \|_{\infty}
& \le & \tilde \mu_{p-j-1}. \label{restterm}
\eeqn

Substituting (\ref{topterm}) and (\ref{restterm}) into (\ref{fullser}) we see that
\beqns
 \| M_{1,p} c_m \|_{\infty} & \le  & \sum_{j=0}^{\ell+1} \tilde \mu_{p-j-1}  +\mu_{p-\ell-2}+(p-\ell-2)A^{p-\ell-3}\epsilon.
\eeqns

\end{proof}

To close we present, in Table~\ref{sigma}, the values of $\tilde \sigma_{p,\ell}$ for $p \ge \ell+5$ where, as before, the computations are based upon the values of  $\mu_{p,1}$ from Table~\ref{numancom}. We note these values give an indication of  the accuracy of multilevel algorithm, from the $5^{th}$ iteration onwards, when applied to a cosine whose frequency falls in the range  $\{2^{\ell},2^{\ell}+1,\ldots,2^{\ell+1}\}.$

\begin{table}[h]
\centering
\label{sigma}
\begin{tabular}{|c|c|c|c|c|c|c|c|c|c|c|} \hline
$p$ & 6 & 7 & 8 & 9 & 10 & 11 & 12 \\ \hline
$\tilde \sigma_{p,1}$   & 6.1(-1) & 4.6(-2) &  1.1(-3) & 1.1(-5)& 1.2(-7) & 3.9(-9) & 1.5(-10) \\ \hline
$\tilde \sigma_{p,2}$   & - & 6.1(-1) & 4.6(-2) &  1.1(-3) & 1.1(-5)& 1.2(-7) & 3.9(-9) \\ \hline
$\tilde \sigma_{p,3}$   & - & -& 6.1(-1) & 4.6(-2) &  1.1(-3) & 1.1(-5)& 1.2(-7) \\ \hline
\end{tabular}
\caption{Values of $\tilde \sigma_{p,\ell}$ for $p\ge \ell+5$  based on the values of $\mu_{p,1}$ in Table~\ref{numancom}}
\end{table}

\section{Numerical experiments} \label{numerics}

In this section we look at four numerical examples. To begin with we examine   $f_{1}=c_1$ and $f_{2}=c_9$  so that we can observe the algorithm treating frequencies similarly once  the sample rate passes the cosine frequency; for $c_{1}$ this occurs after the first level and for $c_{9}$ this occurs after the $4^{th}$ level. Recall that in  Section \ref{fullhalfalg} we derived the sequence $\mu_{p,1}$  whose values are upper bounds on the truncation error $\| T_{\frac{1}{2^{\ell}},p} c_m\|_\infty$ in the case where $m <2^{\ell-1}.$ We note that the leading term in the formula defining $\mu_{p,1}$ is the coefficient $\Gamma_p$ (\ref{gamp}). In Table~\ref{cosone} we compare the convergence of the algorithm applied to $c_{1}$ to the coefficients of both $\Gamma_{p-1}$ and $\mu_{p-1,1}.$ Here we see that the decay rate for $f=c_1$ is almost identically that of $\Gamma_{p-1}$, up until $p=8.$ Beyond $p=8$ the decay rate of the error is influenced by  the coefficients of  cosines, other than $c_1,$ that appear in the error expansion, however we observe that these deviations are well accounted for by the values of $\mu_{p,1}$ and thus, in later stages, these bounds are more reliable.
\begin{table}[h]
\centering
\label{cosone}
\begin{tabular}{|c|c|c|c|c|c|c|c|c|c|c|c|c|} \hline
$p$ & 1 & 2 & 3 & 4 & 5 & 6& 7 & 8 & 9 & 10 & 11 \\ \hline
$\| M_{1,p} c_1\|_\infty$  & 2.0 & 0.99&  0.7 & 0.19& 1.4(-2) & 2.6(-4) & 1.2(-6) & 2.8(-9)& 1.0(-10)& 2.9(-12) & 5.6(-14) \\ \hline
$\Gamma_{p-1}$ & 1 & 0.99&  0.7 & 0.19& 1.4(-2) & 2.6(-4) & 1.3(-6) & 2.5(-9)& 4.6(-13)& 3.5(-17) & 6.5(-22) \\ \hline
$\mu_{p-1,1}$ & - & 1&  0.7 & 0.19& 1.5(-2) & 3.7(-4) & 3.7(-6) & 3.9(-8)& 1.3(-9)& 4.8(-11) & 1.2(-12) \\ \hline
\end{tabular}
\caption{Comparison of multilevel algorithm on $c_1$  and the sequences $\Gamma_{p-1}$ and $\mu_{p-1,1}$ }
\end{table}

In Table~\ref{cosnine} we compare the convergence of the algorithm applied to $c_{9}$ to the coefficients  $\Gamma_{p-4}$ here we observe that the decay of $\| M_{1,p} c_9 \|_\infty$ beyond $4$ levels of iteration consistently tracks $\Gamma_{p-4}.$ This observation reflects the analysis of the previous sections, where convergence happens when the sample rate is greater than the cosine frequency; in the case of $c_{9}$ this occurs at level $4.$

 \begin{table}[h]
\centering
\label{cosnine}
\begin{tabular}{|c|c|c|c|c|c|c|c|c|c|c|c|c|} \hline
$p$ & 1 & 2 & 3 & 4 & 5 & 6& 7 & 8 & 9 & 10 & 11 \\ \hline
$\| M_{1,p} c_9\|_\infty$  & 2.0 & 1.0 &  1.3 & 1.8& 1.0 & 8.0(-1) & 2.6(-1) & 2.4(-2)& 5.7 (-4)& 3.5 (-6) & 7.8 (-9) \\ \hline
$\Gamma_{p-4}$&-&-&-& 1 & 0.99&  0.7 & 0.19& 1.4(-2) & 2.6(-4) & 1.3(-6) & 2.5(-9) \\ \hline
\end{tabular}
\caption{Comparison of multilevel algorithm on $c_9$  and the sequence $\Gamma_{p-4}.$  }
\end{table}

 The third example we consider is the smooth function $f_{3}=\exp(c_1),$ as this allows us to observe the scheme on a function with a full cosine expansion which, appealing to Equation (9.6.19) of \cite{AS}, is given by
\beqns
f(x)=\exp(\cos(2\pi x)) & = & \sum_{k=0}^\infty \widehat{f}(k) \cos(2 \pi k x),\quad {\rm{where}}\,\,\,\widehat{f}(k)  =  \left \{ \begin{array}{ll}
2I_k(1), & k \neq 0, \\
I_0(1), & k=0,
\end{array} \right.
\eeqns
where $I_k$ is the modified Bessel function of order $k$.
In Table~\ref{bessel} we evaluate these Fourier coefficients, which we see decay factorially, and so we may consider  $f_{3}=\exp(c_1)$ as being numerically band-limited.
\begin{table}[h]
\centering
\label{bessel}
\begin{tabular}{|c|c|c|c|c|c|c|c|c|c|c|c|} \hline
$k$ & 0 & 1 & 2 & 3 & 4 & 5& 6 & 7 & 8 & 9 & 10 \\ \hline
$I_k(1)$ & 1.3  & 5.7(-1) & 1.4(-1)& 2.2(-2) & 2.7(-3 & 2.7(-4)& 2.2(-5)& 1.6(-6)& 1.0(-7) & 5.5(-9)& 2.8(-10)\\ \hline
\end{tabular}
\caption{Values of $I_k(1) \ge 10^{-10}$.}
\end{table}

The convergence of the multilevel algorithm applied to $f_{3}$ is illustrated in Table~\ref{bandlimerr}.
 \begin{table}[h]
\centering
\label{bandlimerr}
\begin{tabular}{|c|c|c|c|c|c|c|c|c|c|c|c|} \hline
$p$ & 1 & 2 & 3 & 4 & 5& 6 & 7 & 8 & 9 & 10 & 11 \\ \hline
$\| M_{1,p} f_{3} \|_\infty$ & 2.35  & 1.15 & 1.08& 4.4(-1) & 9.1(-2) & 8.5(-3)& 3.1(-4)& 4.0(-6)& 1.6(-8) & 7.3(-11)& 2.7(-12)\\ \hline
\end{tabular}
\caption{The multilevel method applied to  the numerically band-limited function $f_{3}(x)=\exp(\cos(2 \pi x))$.}
\end{table}

 We now examine how to  develop an error expression, by focusing on the lower frequency cosines, that attempts to closely track the numerical results. We begin by using the analysis of Section~\ref{fullalg} which tells us that for $M_{1,1}c_0 \approx 0$ and for $k \neq 0$, $M_{1,1}c_k \approx c_k-c_0$. Thus the level one error  is
\beqns
\sum_{k=0}^\infty \widehat{f}(k) M_{1,1} c_k & \approx & \sum_{k=1}^\infty \widehat{f}(k) (c_k-c_0) \\
& = & f(x)-\widehat{f}(0)-(f(0)-\widehat{f}(0)) \\
& = & f(x)-f(0).
\eeqns
This is maximised in size when $x=1/2$, with value $e^1-e^{1/2} \approx 2.35$ which is what we observe numerically in Table~\ref{bandlimerr}.\\

For the second iteration we have
\[M_{1,2}c_0 \approx 0,\,\,\, M_{1,2}c_1 \approx \Gamma_1 c_1,\,\,\, {\rm{and}}\,\,\,M_{1,2}c_k \approx c_k-\hpsi(1/2)  c_1,\,\,{\rm{for}} \,\,\,k \ge 2.\]
 Thus, the approximate error at the
second step is
\beqns
\left \| \sum_{k=0}^\infty \widehat{f}(k) M_{1,2} c_k \right \|_{\infty} & \approx & \left \| \widehat{f}(1)\Gamma_1c_1  + \sum_{k=2}^\infty \widehat{f}(k) (c_k-\hpsi(1/2)c_1) \right \|_{\infty} \\
& \le & \Gamma_1 \widehat{f}(1) + (1-\hpsi(1/2)) \sum_{k=2}^\infty \widehat{f}(k) \\
& \approx & 0.99 \times 1.14+0.99 \times(e-1.3-1.14) = 1.4,
\eeqns
which is a reasonable estimate of the error 1.15 in Table~\ref{bandlimerr}. We observe that actually $\Gamma_1 \widehat{f}(1) =0.99  \times 1.14=1.12$ is a much closer approximation, but not an upper bound.

We finish by looking at the third level of iteration, which is more complicated. For the lower frequencies we have
\[
M_{1,3}c_0 \approx 0,\,\,\, M_{1,3}c_1 \approx \Gamma_2 c_1,\,\,\, M_{1,3}c_2 \approx \Gamma_1 c_2\,\,\, {\rm{and}}\,\,\,M_{1,3}c_3 \approx c_{3}-
\hpsi(1/2) \gamma_{2}c_{1}.
\]
For the remaining frequencies $k\ge 4,$ we have the following estimates
\[
M_{1,3}c_k \approx
\begin{cases} c_k - \hpsi(1/2) \gamma_2 c_1 \,\, & \,\,
\textrm{for} \,\,\, k\,\,\,{\rm{odd}}; \\ c_k - c_0 \,\, & \,\,
\textrm{for} \,\,\, k\,\,\,{\rm{divisible}}\,\,\,{\rm{by}}\,\,\,4; \\
c_k - \hpsi(1/2) c_2 \,\, & \,\,
\textrm{otherwise}.
\end{cases}
\]
The calculation we do as above leads us to the following approximate error estimate
\beqns
\left \| \sum_{k=0}^\infty \widehat{f}(k) M_{1,3} c_k \right \|_{\infty} & \approx & \widehat{f}(1)\Gamma_2 + \widehat{f}(2)\Gamma_1 + 2 \sum_{k=3}^\infty \widehat{f}(k) \\
& \approx & 0.7 \times 1.14+0.99 \times 0.28+2*(e-1.3-1.1-0.2) = 1.08,
\eeqns
which is a very accurate approximation of the bound. If we continue this pattern of analysis we can see, from the results displayed in Table \ref{bandlim}, that the expression
\beqn
E_p(f):=\widehat{f}(1)\Gamma_{p-1} + \widehat{f}(2) \Gamma_{p-2} +  \widehat{f}(3) \Gamma_{p-3} + \widehat{f}(4) \Gamma_{p-4}& \approx &
\left \| \sum_{k=0}^\infty \widehat{f}(k) M_{1,p} c_k \right \|_{\infty}   \label{bandlimest}
\eeqn
works as a reasonable estimate to the observed numerical errors.
\begin{table}[h]
\centering
\label{bandlim}
\begin{tabular}{|c|c|c|c|c|c|c|c|c|c|c|c|} \hline
$p$ & 1 & 2 & 3 & 4 & 5& 6 & 7 & 8 & 9 & 10 & 11 \\ \hline
$E_p(f_{3})$ & 2.35  & 1.15 & 1.08 & 4.4(-1) & 9.7(-2) & 1.2(-2)& 6.8(-4)& 1.2(-5)& 5.3(-8) & 1.5(-10) & 5.2(-12) \\
\hline
$\| M_{1,p} f_{3} \|_\infty$ & 2.35  & 1.15 & 1.08& 4.4(-1) & 9.1(-2) & 8.5(-3)& 3.1(-4)& 4.0(-6)& 1.6(-8) & 7.3(-11)& 2.7(-12)\\ \hline
\end{tabular}
\caption{Comparison of actual errors and error estimate based on errors for low frequencies for the numerically band-limited function $f_{3}(x)=\exp(\cos(2 \pi x))$.}
\end{table}

We conclude this section with a second band-limited example
\beqn
f_{4}(x)=\sum_{z \in \Z} \psi(2x-z). \label{pergauss}
\eeqn
This has Fourier coefficients
$$
\widehat{f}(k) = \exp(-2 \pi^2 k^2/4).
$$
If we use (\ref{bandlimest}) again for $p \ge 4$ we obtain the following estimates for the approximation errors $\| M_{1,p} f_{4} \|_\infty$, which we tabulate below for comparison.
\begin{table}[h]
\centering
\label{bandlim2}
\begin{tabular}{|c|c|c|c|c|c|c|c|c|c|c|c|} \hline
$p$ & 1 & 2 & 3 & 4 & 5& 6 & 7 & 8 & 9 & 10 & 11 \\ \hline
$E_p(f_{4})$ & 1.4(-2) & 7.1(-3) & 5.0(-3) & 1.3(-3) & 9.9(-5) & 1.9(-6)& 8.6(-9)& 2.0(-11)& 7.2(-13) & 2.1(-14) & 4.1(-16) \\
\hline
$\| M_{1,p} f_{4} \|_\infty$ & 1.3(-2)  & 7.1(-3) & 5.0(-3) & 1.3(-3) & 9.9(-5) & 1.9(-6)& 9.4(-9)& 2.1(-11)& 7.0(-13) & 2.1(-14)& 5.9(-16)\\ \hline
\end{tabular}
\caption{Comparison of actual errors and error estimate based on errors for low frequencies for the numerically band-limited function $f_{4}$ (\ref{pergauss}).}
\end{table}




\end{document}